\documentclass[a4paper,12pt]{article}
\begin{document}

\title{Integration of vector hydrodynamical partial differential
equations over octonions.}
\author{Ludkovsky S.V.}
\date{29 December 2010}
\maketitle

\begin{abstract}
New technique of integration of certain types of partial
differential equations is developed. For this purpose
non-commutative integration over Cayley-Dickson algebras is used.
Applications to non-linear vector partial differential equations of
Korteweg-de-Vries and Kadomtzev-Petviashvili types and describing
non-isothermal flows of incompressible Newtonian liquids are given.
\end{abstract}

\footnote{key words and phrases: hypercomplex, octonions, non-linear
partial differential equation, non-commutative integration }

\section{Introduction.}
This article is devoted to a method of non-commutative integration
of systems of partial differential equations or a partial
differential equation written in a vector form over Cayley-Dickson
algebras.  This technique is applied to non-linear generalized
Korteweg-de-Vries and Kadomtzev-Petviashvili  and also describing
non-isothermal flows of incompressible Newtonian liquids partial
differential equations written in vector forms with large number of
variables. Dirac had used complexified quaternions to solve
Klein-Gordon's hyperbolic partial differential equation, which is
used in quantum mechanics. This approach is generalized in this
paper for non-linear partial differential equations over the
Cayley-Dickson algebras.
\par The method based on $\bf R$-linear integral
equations was previously used for non-linear partial differential
equations with a real time variable $t$ and real space variables $x,
y$ \cite{ablsigb,polzayjurb}. Using the non-commutative line
integration over the Cayley-Dickson algebras it is spread in this
paper for variables $x, y$ in the Euclidean space ${\bf R}^n$,
$~n\in {\bf N}$. \par Korteweg-de-Vries equations describe a motion
of shallow waters, for example, in channels. Kadomtzev-Petviashvili
and modified Korteweg-de-Vries equations are used in modeling of
physical processes such as non-stationary spread of waves in a
material with dispersion, magneto-hydrodynamical waves in a
non-collision plasma, heat conductivity of a lattice of anharmonic
oscillators.  Moreover, partial differential equations describing
stationary non-isothermal flows of incompressible Newtonian liquids
are very important as well. This method is also applicable to a
resonance interactions in a non-linear material, for example, in the
non-linear optics.
\par In this article previous results of the author on functions of
Cayley-Dickson variable and non-commutative line integrals over
Cayley-Dickson algebras are used
\cite{ludancdnb,ludfov,ludoyst,lulipadecdla}.
\par Henceforward, the notations of previous papers
\cite{ludfov,lulipadecdla} and the book \cite{ludancdnb} are used.

\section{Integration of partial differential equations over Cayley-Dickson algebras.}
\par {\bf 1.1. Notations and Definitions.}
Let ${\cal A}_r$ denote the real Cayley-Dickson algebra with
generators $i_0,...,i_{2^r-1}$ such that $i_0=1$, $~i_j^2=-1$ for
each $j\ge 1$, $~i_ji_k=-i_ki_j$ for each $j\ne k\ge 1$. It is
supposed further, that a domain $U$ in ${\cal A}_r^m$ has the
property that
\par $(D1)$ each projection ${\bf p}_j(U)=:U_j$ is
$(2^r-1)$-connected;
\par $(D2)$ $\pi _{{\sf s},{\sf p},{\sf t}}(U_j)$ is simply
connected in $\bf C$ for each $k=0,1,...,2^{r-1}$, ${\sf s}=i_{2k}$,
${\sf p}=i_{2k+1}$, ${\sf t}\in {\cal A}_{r,{\sf s},{\sf p}}$ and
${\sf u}\in {\bf C}_{{\sf
s},{\sf p}}$, for which there exists $z={\sf u}+{\sf t}\in U_j$, \\
where $e_j = (0,...,0,1,0,...,0)\in {\cal A}_r^m$ is the vector with
$1$ on the $j$-th place, ${\bf p}_j(z) = \mbox{ }^jz$ for each $z\in
{\cal A}_r^m$, $z=\sum_{j=1}^m\mbox{ }^jz e_j$, $\mbox{ }^jz\in
{\cal A}_r$ for each $j=1,...,m$, $m\in {\bf N} := \{ 1,2,3,... \}
$,  $~\pi _{{\sf s},{\sf p},{\sf t}}(V):= \{ {\sf u}: z\in V,
z=\sum_{{\sf v}\in \bf b}w_{\sf v}{\sf v},$ ${\sf u}=w_{\sf s}{\sf
s}+w_{\sf p}{\sf p} \} $ for a domain $V$ in ${\cal A}_r$ for each
${\sf s}\ne {\sf p}\in \bf b$, where ${\sf t}:=\sum_{{\sf v}\in {\bf
b}\setminus \{ {\sf s}, {\sf p} \} } w_{\sf v}{\sf v} \in {\cal
A}_{r,{\sf s},{\sf p}}:= \{ z\in {\cal A}_r:$ $z=\sum_{{\sf v}\in
\bf b} w_{\sf v}{\sf v},$ $w_{\sf s}=w_{\sf p}=0 ,$ $w_{\sf v}\in
\bf R$ $\forall {\sf v}\in {\bf b} \} $, where ${\bf b} := \{
i_0,i_1,...,i_{2^r-1} \} $ is the family of standard generators of
the Cayley-Dickson algebra ${\cal A}_r$. Frequently we take $m=1$.
Henceforward, we consider a domain $U$ satisfying Conditions
$(D1,D2)$ if any other is not outlined. \par The Cayley-Dickson
algebra ${\cal A}_{r+1}$ is formed from the algebra ${\cal A}_r$
with the help of the doubling procedure by generator $i_{2^r}$, in
particular, ${\cal A}_0:=\bf R$ is the real field, ${\cal A}_1=\bf
C$ coincides with the field of complex numbers, ${\cal A}_2=\bf H$
is the skew field of quaternions, ${\cal A}_3$ is the algebra of
octonions, ${\cal A}_4$ is the algebra of sedenions. The skew field
of quaternions is associative, and the algebra of octonions is
alternative. The multiplication of arbitrary octonions $\xi , \eta $
satisfies equations $(1,2)$ below:
\par $(1)$ $(\xi \eta )\eta =\xi (\eta \eta )$,
\par $(2)$ $\xi (\xi \eta )=(\xi \xi )\eta $, \\
that forms the alternative system. The algebra ${\cal A}_r$ is power
associative, that is, $z^{n+m}=z^nz^m$ for each $n, m \in \bf N$ and
$z\in {\cal A}_r$, it is non-associative and non-alternative for
each $r\ge 4$.
\par The Cayley-Dickson algebras are $*$-algebras, that is, there is a
real-linear mapping ${\cal A}_r\ni a\mapsto a^*\in {\cal A}_r$ such
that
\par $(2)$ $a^{**}=a$,
\par $(3)$ $(ab)^*=b^*a^*$ for each $a, b\in {\cal A}_r$.
Then they are nicely normed, that is,
\par $(4)$ $a+a^*=:2 Re (a) \in \bf R$ and
\par $(5)$ $aa^*=a^*a>0$ for each $0\ne a\in {\cal A}_r$.
The norm in it is defined by the equation:
\par $(6)$ $|a|^2:=aa^*$.
\par We also denote $a^*$ by $\tilde a$. Each non-zero Cayley-Dickson
number $0\ne a\in {\cal A}_r$ has the multiplicative inverse given
by $a^{-1}=a^*/|a|^2$.
\par The doubling procedure is as follows. Each $z\in {\cal A}_{r+1}$
is written in the form $z=a+b{\bf l}$, where ${\bf l}^2=-1$, ${\bf
l}\notin {\cal A}_r$, $a, b \in {\cal A}_r$. The addition is
componentwise. The conjugate of a Cayley-Dickson number $z$ is
prescribed by the formula:
\par $(7)$ $z^*:=a^* -b{\bf l}$. \\
The multiplication is given by Equation \par $(8)$ $(\alpha +\beta
{\bf l})(\gamma +\delta {\bf l})=(\alpha \gamma
-{\tilde {\delta }}\beta )+(\delta \alpha +\beta {\tilde {\gamma }}){\bf l}$ \\
for each $\alpha $, $\beta $, $\gamma $, $~ \delta \in {\cal A}_r$,
$~ \xi :=\alpha +\beta {\bf l}\in {\cal A}_{r+1}$, $~ \eta :=\gamma
+\delta {\bf l}\in {\cal A}_{r+1}$.
\par The basis of ${\cal A}_{r+1}$ over $\bf R$ is denoted by
${\bf b}_{r+1} := {\bf b}:= \{ 1, i_1,...,i_{2^{r+1}-1} \} $, where
$i_s^2=-1$ for each $1\le s \le 2^{r+1}-1$, $i_{2^r}:={\bf l}$ is
the additional element of the doubling procedure of ${\cal A}_{r+1}$
from ${\cal A}_r$, choose $i_{2^r+m}= i_m{\bf l}$ for each
$m=1,...,2^r-1$, $i_0:=1$.

\par {\bf 1.2. Operators.} An ${\bf R}$ linear space $X$ which is also
left and right ${\cal A}_r$ module will be called an ${\cal A}_r$
vector space. We present $X$ as the direct sum $X=X_0i_0\oplus ...
\oplus X_{2^r-1} i_{2^r-1}$, where $X_0$,...,$X_{2^r-1}$ are
pairwise isomorphic real linear spaces. Certainly, for $r=2$ this
module is associative: $(xa)b=x(ab)$ and  $(ab)x=a(bx)$ for all
$x\in X$ and $a,b \in {\bf H}$, since the quaternion skew field
${\cal A}_2={\bf H}$ is associative. This module is alternative for
$r=3$: $(xa)a=x(a^2)$ and $(a^2)x=a(ax)$ for all $x\in X$ and $a\in
{\bf O}$, since the octonion algebra ${\bf O}={\cal A}_3$ is
alternative. \par  Let $X$ and $Y$ be two $\bf R$ linear normed
spaces which are also left and right ${\cal A}_r$ modules, where
$1\le r$, such that  $0\le \| ax \|_X = |a| \|x \|_X = \| x a \|_X$
and $\| x+y \|_X \le \| x \|_X + \| y \|_X $ for all $x, y\in X$ and
$a\in {\cal A}_r$.  Such spaces $X$ and $Y$ will be called ${\cal
A}_r$ normed spaces.
\par We say that an ${\cal A}_r$ vector space $X$ is supplied with an ${\cal
A}_r$ valued scalar product, if \par $(f,g) = \sum_{j,k=0}^{2^r-1}
(f_j,g_k)i_j^*i_k$, \\ where $f=f_0i_0+...+f_{2^r-1}i_{2^r-1}$, $
~f, g\in X$, $ ~ f_j, g_j \in X_j$, each $X_j$ is a real linear
space with a real valued scalar product, $(X_j, (*,*))$ is real
linear isomorphic with $(X_k, (*,*))$ and $(f_j,g_k)\in {\bf R}$ for
each $j, k$. The scalar product induces the norm: $\| f \| :=
\sqrt{(f,f)}$.
\par An ${\cal A}_r$ normed space or an
${\cal A}_r$ vector space with ${\cal A}_r$ scalar product complete
relative to its norm will be called an ${\cal A}_r$ Banach space or
an ${\cal A}_r$ Hilbert space respectively.
\par We put $X^{\otimes k} := X\otimes _{\bf R} ... \otimes
_{\bf R} X$ to be the $k$ times ordered tensor product over $\bf R$
of $X$.  By $L_{q,k}(X^{\otimes k},Y)$ we denote a family of all
continuous $k$ times $\bf R$ poly-linear and ${\cal A}_r$ additive
operators from $X^{\otimes k}$ into $Y$. If $X$ and $Y$ are normed
${\cal A}_r$ spaces and $Y$ is complete relative to its norm, then
$L_{q,k}(X^{\otimes k},Y)$ is also a normed $\bf R$ linear and left
and right ${\cal A}_r$ module complete relative to its norm. In
particular, $L_{q,1}(X,Y)$ is denoted also by $L_q(X,Y)$.
\par  If $A\in L_q(X,Y)$ and $A(xb)=(Ax)b$
or $A(bx)=b(Ax)$ for each $x\in X_0$ and $b\in {\cal A}_r$, then an
operator $A$ we call right or left ${\cal A}_r$-linear respectively.
\par An $\bf R$ linear space of left (or right) $k$ times ${\cal A}_r$
poly-linear operators is denoted by $L_{l,k}(X^{\otimes k},Y)$ (or
$L_{r,k}(X^{\otimes k},Y)$ respectively).
\par An $\bf R$-linear operator $A: X\to X$ will be called right or
left strongly ${\cal A}_r$ linear if
\par $A(xb)=(Ax)b$ or
\par $A(bx)=b (Ax)$ for each $x\in X$ and $b\in {\cal A}_r$ correspondingly.
\par {\bf 1.2.1. Examples.} Each right ${\cal A}_r$ linear operator $A: X\to X$ such
that $Ax\in X_0$ for each $x\in X_0$ is strongly right ${\cal A}_r$
linear, since \par $A(xb)= \sum_{j,k=0}^{2^r-1} A(x_ji_jb_ki_k)=
\sum_{j,k=0}^{2^r-1} (Ax_j)i_ji_kb_k$\par $ = \sum_{j,k=0}^{2^r-1}
(Ax_ji_j)i_kb_k =(Ax)b,$ \\ since $Ax_j\in X_j$ up to an ${\bf R}$
linear isomorphism  for each $x=x_0i_0+...+x_{2^r-1}i_{2^r-1}\in X$
and $b=b_0i_0+...+b_{2^r-1}i_{2^r-1}\in {\cal A}_r$ with $x_j\in
X_j$ and $b_j\in {\bf R}$ for each $j=0,...,2^r-1$.
\par Particularly, if $A\in L_r(X,Y)$ is right ${\cal A}_2$ linear with $X$ and $Y$ over ${\cal
A}_2$, i.e. for $X$ and $Y$ over the quaternion skew field ${\bf
H}$, then $A$ is strongly right ${\bf H}$ linear, since the
quaternion skew field is associative and $x(ab) =(xa)b$ for each
$x\in X$ and $a, b \in {\bf H}$.

\par {\bf 1.2.2. Lemma.} {\it If an invertible operator $A$ is either
right ${\cal A}_r$ linear or right strongly ${\cal A}_r$ linear,
then $A^{-1}$ is such also.} \par {\bf Proof.} This follows from the
equalities:
\par $(A^{-1}y)b=(A^{-1}Ax)b=xb =
A^{-1}(A(xb))= A^{-1}((Ax)b)  = A^{-1}(yb)$, \\ where $Ax=y$, either
$x\in X_0$ or $x\in X$ correspondingly and $b\in {\cal A}_r$, since
$A(xb)=(Ax)b$ for each $x\in X_0$ or $x\in X$ respectively and $b\in
{\cal A}_r$.

\par {\bf 1.3. First order partial differential operators.}
We consider an arbitrary first order partial
differential operator $\sigma $ given by the formula
\par $(1)$ $\sigma f= \sum_{j=0}^{2^r-1}
i_j^* (\partial f/\partial z_{\xi (j)}) {\psi }_j$,
\\ where $f$ is a differentiable ${\cal A}_r$-valued function on the domain $U$ satisfying
Conditions 1.1$(D1,D2)$, $~ 2\le r$, $ ~ i_0,...,i_{2^r-1}$ are the
standard generators of the Cayley-Dickson algebra ${\cal A}_r$, $~
\psi _j$ are real constants so that $\sum_j\psi _j^2> 0$, $ ~ \xi :
\{ 0,1,...,2^r-1 \} \to \{ 0,1,...,2^r-1 \}$ is a surjective
bijective mapping, i.e. $\xi $ belongs to the symmetric group
$S_{2^r}$ (see also \S 2 in \cite{lulipadecdla}).
\par Making the substitution $f\mapsto fi_{2^r}$ and using the
embedding ${\cal A}_r\hookrightarrow {\cal A}_{r+1}$ it is always
possible to relate the case of the operator $\sigma $ with $\psi
_0\ne 0$ and $\sigma $ with $\psi _j\ne 0$ only for some $j>0$,
where $i_{2^r}$ is the doubling generator.
\par For an ordered product $\{ \mbox{}_1f...\mbox{}_kf \}_{q(k)}$
of differentiable functions $\mbox{}_sf$ we put
\par $(2)$ $\mbox{}^s\sigma \{ \mbox{}_1f...\mbox{}_kf \}_{q(k)}
= \sum_{j=0}^n  i_j^* \{ \mbox{}_1f...(\partial \mbox{}_sf/\partial
z_{\xi (j)})...\mbox{}_kf \} _{q(k)} \psi _j,$ \\
where a vector $q(k)$ indicates on an order of the multiplication in
the curled brackets (see also \S  2 \cite{ludoyst,ludfov}), so that
\par $(3)$ $\sigma \{ \mbox{}_1f...\mbox{}_kf \}_{q(k)}=\sum_{s=1}^k
\mbox{}^s\sigma \{ \mbox{}_1f...\mbox{}_kf \}_{q(k)}$. \par
\par Symmetrically other operators
\par $(4)$ ${\hat \sigma }f= \sum_{j=0}^{2^r-1}
(\partial f/\partial z_{\xi (j)}) i_j {\psi }_j$, \\
are defined. Therefore, these operators are related:
\par $(5)$ $(\sigma f)^*= {\hat \sigma } (f^*) $.
\par Operators $\sigma $ given by $(1)$ are right ${\cal A}_r$
linear.

\par {\bf 2. Integral operators.} We consider integral operators of the form:
$$(1)\quad {\sf K}(x,y) = {\sf F}(x,y) +\mbox{}_{\sigma } \int_x^{\infty }
{\sf F}(z,y){\sf N}(x,z,y)dz ,$$ where $\sigma $ is an $\bf
R$-linear partial differential operator as in \S 1.3 and
$\mbox{}_{\sigma } \int $ is the non-commutative line integral
(anti-derivative operator) over the Cayley-Dickson algebra ${\cal
A}_r$ from \S \S 2.5 and 2.23 \cite{lulipadecdla} or 4.2.5
\cite{ludancdnb}, where ${\sf F}$ and ${\sf K}$ are continuous
functions with values in the Cayley-Dickson algebra ${\cal A}_r$ or
more generally in the real algebra $Mat_s({\cal A}_r)$ of $s\times
s$ matrices with entries in ${\cal A}_r$. For definiteness we take
$\mbox{}_{\partial /\partial x_0} \int g(z)dz = \int g(z)dz $ and
$\mbox{}_{\sigma } \int g(z)dz$ calculated with the help of the left
algorithm 2.6 \cite{ludoyst,ludfov} or 1.2.6$(LI)$ \cite{ludancdnb}.
In this case $\mbox{}_{\sigma } \int g(z)dz$ is the ${\cal A}_r$
right linear operator in accordance with Lemma 1.2.2.
\par If $\gamma : [a,b]\to {\cal A}_r$ is a function, then
\par $(i)$ $ ~ ~ V_a^b\gamma := \sup_P  |\gamma (t_{j+1})-\gamma (t_j)|$ \\ is
called the variation of $\gamma $ on the segment $[a,b]\subset {\bf
R}$, where the supremum is taken by all finite partitions $P$ of the
segment $[a,b]$, $P= \{ t_0=a<t_1<...<t_n =b \} $, $n\in {\bf N}$. A
continuous function $\gamma : [a,b]\to {\cal A}_r$ with the finite
variation $V_a^b\gamma <\infty $ is called a rectifiable path.
\par Let a domain $U$ be provided with a foliation
by locally rectifiable paths $\{ \gamma ^{\alpha }: ~ \alpha \in
\Lambda \} $, where $\Lambda $ is a set (see below). We take for
definiteness a canonical closed domain $U$ in ${\hat {\cal A}}_r$
satisfying Conditions 1.1$(D1,D2)$ so that $\infty \in U$, where
${\hat {\cal A}}_r = {\cal A}_r \cup \{ \infty \} $ denotes the
one-point compactification of ${\cal A}_r$, $ ~ 2\le r<\infty $.
\par A path $\gamma : <a,b>\to U$ is called locally rectifiable, if
it is rectifiable on each compact segment $[c,e]\subset <a,b>$,
where $<a,b> =[a,b] := \{ t\in {\bf R}: ~ a\le t \le b \} $ or
$<a,b> =[a,b) := \{ t\in {\bf R}: ~ a\le t < b \} $ or $<a,b> =(a,b]
:= \{ t\in {\bf R}: ~ a< t \le b \} $ or $<a,b> =(a,b) := \{ t\in
{\bf R}: ~ a< t < b \} $.
\par A domain $U$ is called foliated by locally rectifiable paths $\{ \gamma
^{\alpha }: ~ \alpha \in \Lambda \} $ if $\gamma : <a_{\alpha
},b_{\alpha }> \to U$ for each $\alpha $ and it satisfies the
following three conditions: \par $(F1)$ $\bigcup_{\alpha \in \Lambda
} \gamma ^{\alpha }(<a_{\alpha },b_{\alpha }>) =U$ and \par $(F2)$
$\gamma ^{\alpha }(<a_{\alpha },b_{\alpha }>)\cap \gamma ^{\beta
}(<a_{\beta },b_{\beta }>) = \emptyset $ for each $\alpha \ne \beta
\in \Lambda $. \\ Moreover, if the boundary $\partial U = cl
(U)\setminus Int (U)$ of the domain $U$ is non-void then \par $(F3)$
$\partial U = (\bigcup_{\alpha \in \Lambda _1} \gamma ^{\alpha
}(a_{\alpha }))\cup
(\bigcup_{\beta \in \Lambda _2} \gamma ^{\beta }(b_{\beta }))$, \\
where $\Lambda _1 = \{ \alpha \in \Lambda : <a_{\alpha },b_{\beta
}>= [a_{\alpha },b_{\beta }> \} $, $\Lambda _2 = \{ \alpha \in
\Lambda : <a_{\alpha },b_{\beta }>= <a_{\alpha },b_{\beta }] \} $.
For the canonical closed subset $U$ we have $cl (U)=U=cl (Int(U))$,
where $cl (U) $ denotes the closure of $U$ in ${\cal A}_v$ and $Int
(U)$ denotes the interior of $U$ in ${\cal A}_v$. For convenience
one can choose $C^1$ foliation, i.e. each $\gamma ^{\alpha }$ is of
class $C^1$. When $U$ is with non-void boundary we choose a
foliation family such that $\bigcup_{\alpha \in \Lambda } \gamma
(a_{\alpha }) =\partial U_1$, where a set $\partial U_1$ is open in
the boundary $\partial U$ and so that $w|_{\partial U_1}$ would be a
sufficient initial condition to characterize a unique branch of an
anti-derivative $w(x)={\cal I}_{\sigma }f(x) = ~ \mbox{}_{\sigma }
\int_{\mbox{}_0x}^x f(z)dz$, \\ where $\mbox{}_0x\in \partial U_1$,
$x\in U$, $\gamma ^{\alpha }(t_0)=\mbox{}_0x$, $\gamma ^{\alpha }
(t)=x$ for some $\alpha \in \Lambda $, $t_0$ and $t\in <a_{\alpha
},b_{\alpha }>$, \par $~ \mbox{}_{\sigma } \int_{\mbox{}_0x}^x
f(z)dz= ~ \mbox{}_{\sigma } \int_{\gamma ^{\alpha }|_{[t_0,t]}}
f(z)dz$.
\par In accordance with Theorems 2.5 and 2.23 \cite{lulipadecdla} or
4.2.5 and 4.2.23 \cite{ludancdnb} the equality
\par $(2)$ $\sigma _x ~ \mbox{}_{\sigma } \int_{\mbox{}_0x}^x
g(z)dz =g(x)$ \\
is satisfied for a continuous function $g$ on a domain $U$ as in \S
1 and a foliation as above.
\par Particularly in the class of ${\cal A}_r$ holomorphic functions
in the domain satisfying Conditions 1.1$(D1,D2)$ this line integral
depends only on initial and final points due to the homotopy theorem
\cite{ludoyst,ludfov}.
\par  We denote by ${\bf
P}={\bf P}(U)$ the family of all locally rectifiable paths $\gamma :
<a_{\gamma },b_{\gamma }>\to U$ supplied with the family of
pseudo-metrics
\par $(3)$ $\rho ^{a,b,c,d}(\gamma ,\omega ) := |\gamma (a)-\omega (c)|+
\inf_{\phi } V_a^b(\gamma (t) - \omega (\phi (t))$ \\
where the infimum is taken by all diffeomorphisms $\phi : [a,b]\to
[c,d]$ so that $\phi (a)=c$ and $\phi (b)=d$, $~ a<b$, $~ c<d$,
$[a,b]\subset <a_{\gamma }, b_{\gamma }>$, $[c,d]\subset <a_{\omega
}, b_{\omega }>$. We take a foliation such that $\Lambda $ is a
uniform space and  the limit
\par $(4)$ $\lim_{\beta \to \alpha }\rho ^{a,b,a,b}(\gamma ^{\beta },\gamma ^{\alpha
})=0$ \\ is zero for each $[a,b]\subset <a_{\alpha }, b_{\alpha }>$.
\par  For example, we can take $\Lambda = {\bf R}^{n-1}$ for a foliation of the entire
Cayley-Dickson algebra ${\cal A}_r$, where $n=2^r$, $t\in {\bf R}$,
$\alpha \in \Lambda $, so that $\bigcup_{\alpha \in \Lambda } \gamma
^{\alpha }([0,\infty )) = \{ z\in {\cal A}_r: ~ Re((z-y)v^*)\ge 0 \}
$ and $\bigcup_{\alpha \in \Lambda } \gamma ^{\alpha }((-\infty ,0])
= \{ z\in {\cal A}_r: ~ Re((z-y)v^*)\le 0 \} $ are two real
half-spaces, where $v, y\in {\cal A}_r$ are marked Cayley-Dickson
numbers and $v\ne 0$. Particularly, we can choose the foliation such
that $\gamma ^{\alpha } (0)=y+\alpha _1v_1+...+\alpha
_{2^r-1}v_{2^r-1}$ and $\gamma ^{\alpha }(t)=tv_0+\gamma ^{\alpha
}(0)$ for each $t\in {\bf R}$, where $v_0,...,v_{2^r-1}$ are $\bf
R$-linearly independent vectors in ${\cal A}_r$.

\par Therefore, the expression $$(5)\quad \mbox{}_{\sigma } \int_x^{\infty }
g(z)dz := \mbox{}_{\sigma } \int_{\gamma ^{\alpha }|_{[t_x,b_{\alpha
}>}} g(z)dz $$ denotes a non-commutative line integral over ${\cal
A}_r$ along a path $\gamma ^{\alpha }$ so that $\gamma ^{\alpha
}(t_x)=x$ and $\lim_{b\to b_{\alpha }}\gamma ^{\alpha }(t)=\infty $
for an integrable function $g$, where $t_x\in <a_{\alpha },b_{\alpha
}>$, $ ~ \alpha \in \Lambda $, $a_{\alpha }=a_{\gamma ^{\alpha }}$,
$b_{\alpha }=b_{\gamma ^{\alpha }}$. It is sometimes convenient to
use the line integral
$$(6)\quad \mbox{}_{\sigma } \int_x^{-\infty }
g(z)dz := \mbox{}_{\sigma } \int_{\gamma ^{\alpha }|_{<a_{\alpha
},t_x]}} g(z)dz ,$$ when $\lim_{a\to a_{\alpha }}\gamma ^{\alpha
}(t)= \infty $.
\par If $\mbox{}_{\sigma }\int g(z)dz=\sigma ^{-1}g$ corresponds to
the right (or left) algorithm of the non-commutative line
integration over the Cayley-Dickson algebra ${\cal A}_r$, then
$i_{2^r}[(\mbox{}_{\sigma }\int g(z)dz)i^*_{2^r}]$ corresponds to
the left (or right correspondingly) algorithm such that \par
$i_{2^r}[(\sigma f)i^*_{2^r}]=\sum_{j=0}^{2^r-1}(\partial
f^*/\partial z_{\xi (j)})i_j\psi _j={\hat \sigma }(f^*)$, \\ where
$f$ is a differentiable ${\cal A}_r$ valued function.
\par For example, one can take a function
$g(z)$ as ${\sf F}(z,y){\sf K}(x,z)$ depending on additional
parameters $x, y, t,...$.
\par A function ${\sf F}$ has the decomposition
\par $(7)$ ${\sf F} = \sum_{j=0}^{2^r-1} {\sf F}_ji_j$, \\ where each function ${\sf F}_j$
is real-valued for a function ${\sf F}$.
\par Put for convenience $\sigma ^0=I$, where $I$ denotes the unit
operator, $ ~ \sigma ^m$ denotes the $m$-th power of $\sigma $ for
each non-negative integer $0\le m\in {\bf Z}$.

\par {\bf 3. Proposition.} {\it A family ${\cal D}_r$ of all differential operators
with constant ${\cal A}_r$ coefficients is a power associative real
algebra with a center $Z({\cal D}_r)$ consisting of all differential
operators with real coefficients and with a unit element $I$.}
\par {\bf Proof.} Let $A$ be a differential operator with constant
${\cal A}_r$ coefficients, then it can be written in the form:
\par $(1)$ $Af=\sum_{j,s=0}^{2^r-1} [A_{j,s}\pi _sf] i_j^*$, \\ where $A_{j,s}$ is a differential
operator with real coefficients for each $j$:
\par $(2)$ $A_{j,s}=\sum_{\alpha } a_{j,s,\alpha }\partial ^{\alpha }$, \\
$a_{j,s,\alpha }\in {\bf R}$ for each $\alpha = (\alpha
_0,...,\alpha _{2^r-1})$, $~ \partial ^{\alpha }= \frac{\partial
^{|\alpha |}}{\partial x_0^{\alpha _0}...\partial x_{2^r-1}^{\alpha
_{2^r-1}}}$, $~ |\alpha | =\alpha _0+...+\alpha _{2^r-1}$, $~\alpha
_k=0,1,2,...$ is a non-negative integer for each $k=0,...,2^r-1$,
$~\pi _s: X\to X_s$ denotes the ${\bf R}$-linear projection operator
(see \S 1.2),
\par $(3)$ $Af=\sum_{j,s=0}^{2^r-1} (A_{j,s}f_s)i_j^*$ \\
for each $f$ in a domain $Dom (A)$ of the operator $A$, where
\par $(4)$ $f=\sum_{s=0}^{2^r-1}f_si_s$, \\ $f_s$ is a real-valued function for
each $s$, $~f_s=\pi_s f$. If $P_n(z)$ is a polynomial with ${\cal
A}_r$ coefficients of the variable $z\in {\cal A}_r$ (see also \S
2.1 \cite{ludfov,ludoyst}), then
\par $(5)$ $AP_n(z)=\sum_{j,s=0}^{2^r-1}\sum_{\alpha , |\alpha |\le n
} a_{j,s,\alpha }[\partial ^{\alpha } (P_n)_s(z)]i_j^*$,
\\ consequently, $P_n$ belongs to the domain $Dom (A)$ of the operator $A$.
Thus $Dom (A)$ is non-void for each differential operator $A\in
{\cal D}_r$. From Formulas $(1-4)$ we get, that the sum $A+B$ and
the product $BA$ of any differential operators $A, B\in {\cal D}_r$
belongs to ${\cal D}_r$, since
\par $(6)\quad BAf=B(Af)=\sum_{k,j,s=0}^{2^r-1} (-1)^{sign
(j)}(B_{k,j}A_{j,s}f_s)i_k^*=Cf$ with
\par $(7)$ $C_{k,s}=\sum_{j=0}^{2^r-1}(-1)^{sign (j)}B_{k,j}A_{j,s}$, \\
where $sign (x)=0$ for $x=0$, $ ~ sign (x)=1$ for $x>0$, $~
sign(x)=-1$ for $x<0$. Certainly $tA\in{\cal D}_r$ for each $t\in
{\bf R}$ and $A\in {\cal D}_r$. From Formulas 4$(7-9)$
\cite{lulipadecdla} or 4.2.4$(7-9)$ \cite{ludancdnb} we deduce that
\par $(8)$ $A^2 f = A(Af) = \sum_{j=0}^{2^r-1} (A_j^2f) i_j^2
+ \sum_{j=0}^{2^r-1} [(A_0A_j+A_jA_0) f] i_j^* $ \\ for each $f\in Dom
(A^2)$, since $I x = \sum_{j=0}^{2^r-1} \pi _j(x)i_j$ for each $x\in X$,
where \par $(9)$ $A_jf=\sum_{s=0}^{2^r-1}A_{j,s}f_s$, \\
since $(i_j^*)^2=i_j^2$. Therefore, by induction we infer that
$A^{2k}=(A^2)^k$ and $A^{2k+1}=(A^2)^kA=A(A^2)^k$ for each natural
number $k$, where $A^0=I$ denotes the unit operator, i.e. $If=f$ for
each function $f: {\cal A}_r\to {\cal A}_r$. Thus the algebra ${\cal
D}_r$ is power associative: \par $(10)$ $A^kA^m=A^{k+m}$ for each
natural numbers $k, m\in {\bf N}$ and every $A\in {\cal D}_r$.
Moreover it contains the unit element $I$.
\par If $A\in Z({\cal D}_r)$, i.e. an operator $A$ commutes with
every differential operator $B\in {\cal D}_r$, $[A,B]:=AB-BA=0$,
then it commutes with all generators $\{ i_j, \partial _{x_k}i_j^* :
~ j, k =0,...,2^r-1 \} $ of the algebra ${\cal D}_r$. From Formulas
$(1-4,6,7)$ it follows that $a_{j,s,\alpha }=0$ for $j>0$, where
$a_{j,s,\alpha }\in {\bf R}$ for each $j,s$ and $\alpha $. That is,
all the coefficients $a_{j,s,\alpha }i_j^*$ of the operator $A$ are
real, since $i_0=1$.
\par {\bf 3.1. Corollary.} {\it Let $A\in {\cal D}_r$ be a
differentiable operator with constant ${\cal A}_r$ coefficients and
$\pi _0: X\to X_0$ be the ${\bf R}$ linear projection (see \S 1.2),
then $[(I-\pi _0)A]^{2k}$ is with real coefficients for each natural
number $k\in {\bf N}$.}
\par {\bf Proof.}  In view of Formulas 3$(1,9)$ the equality is
valid: \par $(1)$ $(I-\pi _0)A f = \sum_{j=1}^{2^r-1} (A_jf)i_j^*$,
 consequently,
\par $(2)$ $[(I-\pi _0)A]^2 f = \sum_{j=1}^{2^r-1} (A_j^2f) i_j^2$.
\\ That is, $[(I-\pi _0)A]^2$ is the differential operator with real coefficients,
since $i_j^2\in {\bf R}$ for each $j$.
\par {\bf 3.2. Remark.} Mention that in the considered case $\psi
_0=0$ and the operator $\sigma ^2$ is with real coefficients due to
Corollary 3.1, hence the square of the anti-derivative operator
$(\mbox{}_{\sigma }\int ^x)^2$ or $\sigma _x^{-2}$ is also with real
coefficients, when each $\psi _j$ is a real constant and $\psi
_0=0$.

\par {\bf 4. Proposition.} {\it Let $$(1)\quad \lim_{z\to \infty } ~ \mbox{}^1\sigma _z^k ~ \mbox{}^2\sigma
_x^s ~ \mbox{}^2\sigma _z^n{\sf F}(z,y){\sf K}(x,z)=0$$ for each $x,
y$ in a domain $U$ satisfying Conditions 1.1$(D1,D2)$ with $\infty
\in U$ and every non-negative integers $0\le k, s, n\in {\bf Z}$
such that $k+s+n\le m $. Suppose also that $\mbox{}_{\sigma }
\int_x^{\infty } \partial ^{\alpha }_x\partial ^{\beta }_y \partial
^{\omega }_z [{\sf F}(z,y) {\sf K}(x,z)] dz $ converges uniformly by
parameters $x, y $ on each compact subset $W\subset U\subset {\cal
A}_r^2$ for each $|\alpha |+|\beta |+|\omega |\le m$, where $\alpha
=(\alpha _0,...,\alpha _{2^r-1})$, $|\alpha |=\alpha _0+...+\alpha
_{2^r-1}$, $\partial ^{\alpha }_x=\partial ^{|\alpha |}/\partial
x_0^{\alpha _0} ...\partial x_{2^r-1}^{\alpha _{2^r-1}}$. Then the
non-commutative line integral $\mbox{}_{\sigma } \int_x^{\infty }
{\sf F}(z,y){\sf K}(x,z)dz $ from \S 2 satisfies the identities:
$$(2)\quad \sigma ^m_x ~ \mbox{}_{\sigma}\int_x^{\infty } {\sf F}(z,y){\sf K}(x,z)dz =
~ \mbox{}^2\sigma ^m_x ~ \mbox{}_{\sigma}\int_x^{\infty } {\sf
F}(z,y) {\sf K}(x,z)dz + A_m({\sf F},{\sf K})(x,y),$$
$$(3)\quad \mbox{}^1\sigma ^m_z~ \mbox{}_{\sigma}\int_x^{\infty }  {\sf F}(z,y) {\sf K}(x,z)dz =
(-1)^m ~ \mbox{}^2\sigma ^m_z ~\mbox{}_{\sigma}\int_x^{\infty } {\sf
F}(z,y){\sf K}(x,z)dz + B_m({\sf F},{\sf K})(x,y),$$ where
$$(4)\quad A_m({\sf F},{\sf K})(x,y)= - ~\mbox{}^2\sigma
^{m-1}_x [{\sf F}(x,y){\sf K}(x,z)]|_{z=x} + \sigma _x~ A_{m-1}({\sf
F},{\sf K})(x,y)$$ for $m\ge 2$, $$(5) \quad B_m({\sf F},{\sf
K})(x,y) = (-1)^m ~ \mbox{}^2\sigma ^{m-1}_z {\sf F}(x,y){\sf
K}(x,z)|_{z=x} + [\mbox{}^1\sigma _z~ B_{m-1}({\sf F}(z,y),{\sf
K}(x,z))]|_{z=x}$$ for $m\ge 2$,
\par $(6)\quad  ~ A_1({\sf F},{\sf K})(x,y)= B_1({\sf F},{\sf K})(x,y)= - {\sf F}(x,y){\sf K}(x,x)$, \\
 $ ~ \sigma _x$ is an operator $\sigma $ acting by the
variable $x\in U\subset {\cal A}_r$.}
\par {\bf Proof.} From the conditions of this proposition and the
theorem about differentiability of improper integrals by parameters
we have the equality \par $\mbox{}_{\sigma } \int_x^{\infty
}\partial ^{\alpha }_x\partial ^{\beta }_y  \partial ^{\omega }_z
[{\sf F}(z,y) {\sf K}(x,z)] dz = \partial ^{\alpha }_x\partial
^{\beta }_y ~ \mbox{}_{\sigma } \int_x^{\infty }  \partial ^{\omega
}_z [{\sf F}(z,y) {\sf K}(x,z)] dz$ \\ for each $|\alpha |+|\beta
|+|\omega |\le m$ (see also Part IV, Chapter 2, \S 4 in
\cite{kamyn}).
\par In view of  Theorems 2.5 and 2.23 and Corollary 2.6 \cite{lulipadecdla}
or 4.2.5 and 4.2.23 and Corollary 4.2.6 \cite{ludancdnb}  the
equalities
$$(7)\quad \sigma _x ~ \mbox{}_{\sigma} \int_x^{\infty }g(z)dz = -g(x)\mbox{  and}$$
$$(8)\quad \mbox{}_{\sigma} \int_{\mbox{}_0x}^x [\sigma _zf(z)]dz
=f(x)-f(\mbox{}_0x)$$ are satisfied for each continuous function $g$
and continuously differentiable functions $g$ and $f$, where
$\mbox{}_0x$ is a marked point in $U$,
$$(9)\quad \mbox{}^1\sigma _z ~ \mbox{}_{\sigma } \int_x^{\infty }
{\sf F}(z,y){\sf K}(x,z)dz := \sum_{j=0}^{2^r-1} ~ \mbox{}_{\sigma }
\int_x^{\infty } \{ i_j^* [(\partial {\sf F}(z,y)/\partial z_{\xi
(j)}){\sf K}(x,z)]\psi _j \} dz \mbox{  and}$$
$$(10)\quad \mbox{}^2\sigma _z ~ \mbox{}_{\sigma} \int_{\mbox{}_0x}^x {\sf F}(z,y){\sf K}(x,z)dz
:= \sum_{j=0}^{2^r-1} ~ \mbox{}_{\sigma} \int_{\mbox{}_0x}^x \{
i_j^* [ {\sf F}(z,y)(\partial {\sf K}(x,z)/\partial z_{\xi
(j)})]\psi _j \} dz \mbox{ and}$$
$$(11)\quad ~
\mbox{}^2\sigma _x ~ \mbox{}_{\sigma}\int_x^{\infty } {\sf
F}(z,y){\sf K}(x,z)dz := \sum_{j=0}^{2^r-1}
\mbox{}_{\sigma}\int_x^{\infty } \{ i_j^* [{\sf F}(z,y)(\partial
{\sf K}(x,z)/\partial x_{\xi (j)})] \psi _j \} dz.$$ Therefore, from
Equalities $(7,8)$, 1.3$(3)$ and 2$(5)$ and Condition $(1)$ one
gets:
$$(12)\quad \sigma _x~ \mbox{}_{\sigma}\int_x^{\infty }  {\sf F}(z,y) {\sf K}(x,z)dz =
~ \mbox{}^2\sigma _x ~ \mbox{}_{\sigma}\int_x^{\infty } {\sf
F}(z,y){\sf K}(x,z)dz - {\sf F}(x,y){\sf K}(x,x),$$ since ${\sf
F}(z,y){\sf K}(x,z)|_x^{\infty } = - {\sf F}(x,y){\sf K}(x,x)$, that
demonstrates Formula $(2)$ for $m=1$ and $A_1=-{\sf F}(x,y){\sf
K}(x,x)$. The induction by $p=2,...,m$ gives:

$$(13)\quad \sigma ^p_x~ \mbox{}_{\sigma}\int_x^{\infty }  {\sf F}(z,y) {\sf K}(x,z)dz
=$$  $$ ~ \sigma _x ~  [\mbox{}^2\sigma ^{p-1}_x
~\mbox{}_{\sigma}\int_x^{\infty } {\sf F}(z,y){\sf K}(x,z)dz] +
\sigma _x ~ A_{p-1}({\sf F},{\sf K})(x,y)$$
$$= ~ \mbox{}^2\sigma ^p_x ~\mbox{}_{\sigma}\int_x^{\infty }
{\sf F}(z,y){\sf K}(x,z)dz$$ $$ +  ~ \mbox{}^2\sigma
^{p-1}_xA_1({\sf F}(x,y),{\sf K}(x,z)|_{z=x}) + \sigma _x~
A_{p-1}({\sf F},{\sf K})(x,y).$$ Substituting in these expressions
$A_1$ one gets Formulas $(2,4)$ for each $m\ge 2$. Then from
Formulas $(7,8)$ and Condition $(1)$ we infer also that

$$(14)\quad \mbox{}^1\sigma _z~ \mbox{}_{\sigma}\int_x^{\infty }  {\sf F}(z,y) {\sf K}(x,z)dz =
- ~ \mbox{}^2\sigma _z ~ \mbox{}_{\sigma}\int_x^{\infty } {\sf
F}(z,y){\sf K}(x,z)dz +{\sf F}(z,y){\sf K}(x,z)|_x^{\infty }$$ $$ =
-{\sf F}(x,y){\sf K}(x,x) - ~ \mbox{}^2\sigma _z ~
\mbox{}_{\sigma}\int_x^{\infty } {\sf F}(z,y){\sf K}(x,z)dz,$$ that
gives Formulas $(3)$ for $m=1$ and $(6)$. Then we deduce Formulas
$(3,5)$ by induction on $p=2,...,m$:
$$(15)\quad \mbox{}^1\sigma ^p_z~ \mbox{}_{\sigma}\int_x^{\infty }  {\sf F}(z,y) {\sf K}(x,z)dz
=$$  $$(-1)^{p-1} ~\mbox{}^1\sigma _z~ [ \mbox{}^2\sigma ^{p-1}_z
~\mbox{}_{\sigma}\int_x^{\infty } {\sf F}(z,y){\sf K}(x,z)dz ] +
\mbox{}^1\sigma _z~ B_{p-1}({\sf F}(z,y),{\sf K}(x,z))|_{z=x}$$
$$=(-1)^p ~ \mbox{}^2\sigma ^p_z ~\mbox{}_{\sigma}\int_x^{\infty }
{\sf F}(z,y){\sf K}(x,z)dz$$ $$ + (-1)^p ~\mbox{}^2\sigma
^{p-1}_z{\sf F}(z,y){\sf K}(x,z)|_{z=x}+ \mbox{}^1\sigma _z~
B_{p-1}({\sf F}(z,y),{\sf K}(x,z))|_{z=x} .$$
\par {\bf 4.1. Remark.} The center of the Cayley-Dickson algebra ${\cal A}_r$ for $r\ge 2$ is the
real field ${\bf R}$. \par Let $Mat_s({\cal A}_r)$ denote the left
and right ${\cal A}_r$ module of square $s\times s$ matrices with
entries in ${\cal A}_r$. It is possible to consider the quaternion
skew field ${\bf H}_{J,K,L}$ with generators $J, K, L$ realized as
$4\times 4$ square real matrices. Putting $i_jY = (y_{k,l}i_j) =
Yi_j$ for each real matrix $Y$ with elements $y_{k,l}\in {\bf R}$
for each $k, l$ and every generator $i_j$ of the Cayley-Dickson
algebra we naturally obtain the quaternionified algebra $({\cal
A}_r)_{{\bf H}_{J,K,L}}$  (see also \cite{kansol,kurosh}). On the
other hand, ${\bf H}\subset Mat_4({\bf R})$, consequently,
$Mat_s(({\cal A}_r)_{{\bf H}_{J,K,L}})\subset Mat_{4s}({\cal A}_r)$.

\par {\bf 5. Corollary.} {\it If suppositions of Proposition 4 are
satisfied, then
\par $(1)$ $A_2({\sf F},{\sf K})(x,y) = - \sigma _x [{\sf F}(x,y){\sf K}(x,x)] - ~\mbox{}^2\sigma
_x[{\sf F}(x,y){\sf K}(x,z)]|_{z=x}$,
\par $(2)$ $A_3({\sf F},{\sf K})(x,y) = - \sigma ^2_x [{\sf F}(x,y){\sf K}(x,x)]$\par $ - \sigma _x (~\mbox{}^2\sigma
_x [{\sf F}(x,y){\sf K}(x,z)]|_{z=x}) - ~\mbox{}^2\sigma ^2_x[{\sf
F}(x,y){\sf K}(x,z)]|_{z=x}$,
\par $(3)$ $B_2({\sf F},{\sf K})(x,y) = - \mbox{}^1\sigma _x[{\sf F}(x,y){\sf K}(x,x)]
+ ~\mbox{}^2\sigma _z[{\sf F}(x,y){\sf K}(x,z)]|_{z=x}$,
\par $(4)$ $B_3({\sf F},{\sf K})(x,y) = - \mbox{}^1\sigma ^2_x[{\sf F}(x,y){\sf K}(x,x)]
$\par  $ + \mbox{}^1\sigma _x(\mbox{}^2\sigma _z[{\sf F}(x,y){\sf
K}(x,z)]|_{z=x}]) - ~\mbox{}^2\sigma ^2_z[{\sf F}(x,y){\sf
K}(x,z)]|_{z=x} $.
\par $(5)$ $A_2({\sf F},{\sf K})(x,y)-B_2({\sf F},{\sf K})(x,y) = - 2 ~\mbox{}^2\sigma _x [{\sf F}(x,y)
{\sf K}(x,x)]$, \\ where $\sigma _x{\sf K}(x,x)=[\sigma _x{\sf
K}(x,z)+\sigma _z{\sf K}(x,z)]|_{z=x}$,
\par $(6)$ $A_3({\sf F},{\sf K})(x,y)-B_3({\sf F},{\sf K})(x,y)=
- (3 ~\mbox{}^2\sigma ^2_x + ~ \mbox{}^2\sigma _x ~ \mbox{}^2\sigma
_z+ 2~ \mbox{}^2\sigma _z ~ \mbox{}^2\sigma _x) [{\sf F}(x,y){\sf
K}(x,z)]|_{z=x}$\par $ - (2 ~\mbox{}^1\sigma _x ~\mbox{}^2\sigma _x
+ ~\mbox{}^2\sigma _x ~\mbox{}^1\sigma _x) [{\sf F}(x,y){\sf
K}(x,x)]$.
\par Particularly, if either $p$ is even and $\psi _0=0$, or ${\sf F}\in Mat_s({\bf
R})$ and ${\sf K}\in Mat_s({\cal A}_r)$, then \par $(7)$
$\mbox{}^2\sigma _x^p[{\sf F}(z,y){\sf K}(x,z)]= {\sf F}(z,y)\sigma
_x^p{\sf K}(x,z)$ and $\mbox{}^2\sigma _z^p[{\sf F}(z,y){\sf
K}(x,z)]= {\sf F}(z,y)\sigma _z^p{\sf K}(x,z)$.}
\par {\bf Proof.} Formulas $(1-6)$ follow from Equalities 4$(4-6)$.
In particular, if $p$ is even and $\psi _0=0$, $p=2k$, $k\in {\bf
N}$, then in accordance with Corollary 3.1 and Formulas 4$(7-9)$
\cite{lulipadecdla} or 4.2.4$(7-9)$ \cite{ludancdnb} \par $\sigma
_x^pf(x)=A^kf(x)$ \\ for $p$ times differentiable function $f: U\to
{\cal A}_r$, where \par $Af=\sum_jb_j\partial ^2f(x)/\partial
x_j^2$, $b_j=i_{\xi ^{-1}(j)}^2\in {\bf R}$. \par Evidently the
operators $\mbox{}^2\sigma _x^p$ and $\mbox{}^2\sigma _z^p$ commute
with the left multiplication on ${\sf F}(z,y)\in Mat_s({\bf R})$ in
accordance with \S 4.1, that is $\mbox{}^2\sigma _x^p[{\sf
F}(z,y){\sf K}(x,z)]= {\sf F}(z,y)\sigma _x^p{\sf K}(x,z)$ and
$\mbox{}^2\sigma _z^p[{\sf F}(z,y){\sf K}(x,z)]= {\sf F}(z,y)\sigma
_z^p{\sf K}(x,z)$.

\par {\bf 6. Remark.} Expressions of
functions $A_m$ and $B_m$ depend not only on ${\sf F}$ and ${\sf
K}$, but also on $\sigma $, that is on coefficients $\psi _j$ for
$j=0,...,2^r-1$ with $2\le r$ in accordance with \S \S 4 and 5.

\par {\bf 7. Method of non-commutative integration of vector partial
differential equations.} \par Let an equation over the
Cayley-Dickson algebra ${\cal A}_r$ be given in the non-commutative
line integral form:
$$(1)\quad {\sf K}(x,y) = {\sf F}(x,y) + {\sf p} ~ \mbox{}_{\sigma }
\int_x^{\infty } {\sf F}(z,y)N(x,z,y)dz ,$$ where ${\sf K}$, ${\sf
F}$ and ${\sf N}$ are continuous integrable functions as in \S 2,
${\sf p} ~ \in {\bf R}\setminus \{ 0 \} $ is a non-zero real
constant. These functions may depend on additional parameters $t,
\tau,...$. \par  The first step consists in a concretization of a
function ${\sf N}$ and its expression throughout ${\sf F}$ and ${\sf
K}$. Suppose that an operator
\par $(2)$ $(I-{\sf A}_x){\sf F}(x,y)={\sf K}(x,y)$ is invertible,
so that $(I-{\sf A}_x)^{-1}$ is continuous, where $I$ denotes the
unit operator, while \par ${\sf A}_x= - {\sf p} ~ \mbox{}_{\sigma }
\int_x^{\infty } {\sf F}(z,y){\sf N}(x,z,y)dz $ \\ is an operator
acting by variables $x$. \par On the second step two ${\bf
R}$-linear differential or partial differential operators $L_s$ over
the Cayley-Dickson algebra ${\cal A}_r$ are given:
\par $(3)$ $L_sf = \sum_j i_j^* (L_{s,j}f)$, \\ where
$f$ is a differentiable function in the domain of $L_s$, $~L_{s,j}g$
is real-valued function for each $ord (L_s)$ times differentiable
real-valued function $g$ in the domain of $L_{s,j}$ for every $j$
(see also Formulas 3$(1-4)$), where $ord (L_s)$ denotes the order of
the partial differential operator $L_s$. On a function ${\sf F}$ the
conditions either:
\par $(4)$ $L_s{\sf F} =0$ \\ for each $s=1,2$, or
\par $(5)$ $\sum_{j\in \Psi _k} i_j^* [~ \mbox{}_jc(L_{s,0}{\sf F})+L_{s,j}{\sf F}] =0$ \\
for each $s=1,2$ and $1\le k \le m$ are imposed, where $\mbox{}_jc$
are constants $\mbox{}_jc\in {\cal A}_r$, $~\Psi _k\subset \{
0,1,...,2^r-1 \} $ for each $k$, $ ~ \bigcup_k \Psi _k = \{ 0, 1,
..., 2^r-1 \} $, $~ \Psi _k \cap \Psi _l = \emptyset $ for each
$k\ne l$, $~1\le m\le 2^r$.  Coefficients $\mbox{}_jc$ or operators
$L_{s,j}$ may be zero for some $j$. \par On the third step a
function ${\sf K}$ is calculated from Equation $(2)$. \par Acting by
the operator $L_s$ from the left on $(2)$ and using Conditions
either $(4)$ or  $(5)$ one gets either
\par $(6)$ $L_s[(I-A_x){\sf K}]=0$ or
the equalities:
\par $(7)$ $\sum_{j\in \Psi _k} i_j^*
\{ ~ \mbox{}_jc L_{s,0}[(I-{\sf A}_x){\sf K}] + L_{s,j}[(I-{\sf
A}_x){\sf K}]\} =0$ for each $k=1,...,m$.
\\ Therefore, using Conditions either $(4)$ or $(5)$ we infer that
\par $(8)$ $(I-{\sf A}_x)(L_s{\sf K})=R_s({\sf K})$ for $s=1,2$, \\ where an operator
\par $(9)$ $R_s(f)=(I-{\sf A}_x)(L_sf)-L_s[(I-{\sf A}_x)f]$ \\
is formed with the help of commutators $[A,B]=AB-BA$ and
anti-commutators $\{ A,B \} = AB+BA$ of operators $(I-({\sf
A}_x)_0)$, $~ ({\sf A}_x)_k$, $~L_{s,j}$, $~k, j=0,...,2^r-1$. A
function ${\sf N}$ and operators $L_1$ and $L_2$ are chosen such
that
\par $(10)$ $R_s({\sf K})=(I-{\sf A}_x)M_s({\sf K})$ for $s=1, 2$, \\ where $M_s({\sf K})$ are functionals of
${\sf K}$ which generally may be non-${\bf R}$-linear. In view of
Condition $(2)$ the function ${\sf K}$ must satisfy partial
differential equations
\par $(11)$ $L_s{\sf K}-M_s({\sf K})=0$ for $s= 1, 2$, \\ which generally may be non-${\bf
R}$-linear. Thus each solution ${\sf K}$ of the $\bf R$-linear
integral equation $(1)$ is also the solution of partial differential
equations $(11)$. Frequently particular cases are considered, when
Equations $(4)$ correspond to an eigenvalue problem for $s=1$ and to
an evolution in time problem for $s=2$. Generally operators
$L_{s,0},...,L_{s,2^r-1}$ can be chosen ${\cal A}_r$ vector
independent and $\mbox{}_jc\in {\bf R}i_{k(j)}$ for each $j$ and
$s=1,2$ with $k=k(j)\in \{ 0,1,..,2^r-1 \} $.
\par The Euclidean space ${\bf R}^{2^r}$ is the real shadow of the
Cayley-Dickson algebra ${\cal A}_r$, that is, by the definition
${\cal A}_r$ considered as the ${\bf R}$ linear space is isomorphic
with ${\bf R}^{2^r}$. The Lebesgue (non-negative) measure $\mu $ on
the Borel $\sigma $-algebra ${\cal B}({\bf R}^{2^r})$ of the
Euclidean space ${\bf R}^{2^r}$ induces the Lebesgue measure on
${\cal B}({\cal A}_r)$. Therefore, the Hilbert space $X=L^2({\cal
A}_r,\mu ,{\cal A}_r)$  of all Lebesgue measurable functions $f:
{\cal A}_r\to {\cal A}_r$ with integrable square module $|f|^2$,
i.e. \par $(HN)\quad \| f\|^2 := \int_{{\cal A}_r} |f(z)|^2 \mu
(dz)<\infty $, and with the ${\cal A}_r$ valued scalar product \par
$(SP)\quad (f,g) := \int_{{\cal A}_r} f^*(z)g(z) \mu (dz)$ \\
exists. Analogously the Hilbert space $X=L^2({\cal A}_r,\mu
,Mat_s({\cal A}_r))$ is defined with $\sum_{j,k=0}^s f_{j,k}^*(z)
g_{j,k}(z)$ instead of $f^*(z)g(z)$ in the integral in Formula
$(SP)$, where $f_{j,k}\in {\cal A}_r$ denotes a matrix element at
the intersection of row $j$ with column $k$ of a matrix $f\in
Mat_s({\cal A}_r)$.
\par Let $A$ be an ${\bf R}$ linear ${\cal A}_r$ additive operator
$A: {\sf D}(X) \to Y$, where ${\sf D}(A)$ is a domain of $A$ dense
in $X$, ${\sf D}(A)\subset X$, $X$ and $Y$ are Hilbert spaces over
${\cal A}_r$. Then its adjoint operator $A^*$ is defined on a domain
consisting of all those vectors $y\in Y$ such that for some vector
$z\in X$ the equality $(x,z) = (Ax,y)$ is valid for all $x\in {\sf
D}(A)$. For such $y\in Y$ we put $A^*y=z$. If $A^*=A$ we say that
$A$ is self-adjoint. Thus $(Ax,y) = (x,A^*y)$ for all $x\in {\sf
D}(A)$ and $y\in {\sf D}(A^*)$.
\par If an operator $A$ is self-adjoint $A=A^*$, then $(Ax,y)
=(x,Ay)$ for all $x, y\in {\sf D}(A)$.
\par When ${\sf D}(A)$ is dense in $X$ and $(Ax,y) = (x,Ay)$ for all
$x, y\in {\sf D}(A)$, we say that $A$ is symmetric. A self-adjoint
operator is maximal symmetric.
\par We define a self-adjoint operator $A$ to be positive, when
$(Ax,x)\ge 0$ for each $x\in {\sf D}(A)$.
\par An ${\bf R}$ linear ${\cal A}_r$ additive operator $A$ is called
invertible if it is densely defined and one-to-one and has dense
range ${\cal R}(A)$.
\par The operator $A^*A$ is self-adjoint and positive. If $A$ is invertible, then $(A^*A)^{-1}A^*\subseteq
A^{-1}$ (see also \cite{cadringb,lujmsalop}). \par If an expression
of the form
\par $(12)$ $\sum_k [(I-A_x)~\mbox{}_kf(x,y)]~\mbox{}_kg(y)=u(x,y)$ \\ will appear on a domain $U$,
which need to be inverted we consider the case when \par $(RS)$
$(I-A_x)$ is either right strongly ${\cal A}_r$ linear, or right
${\cal A}_r$ linear and $~\mbox{}_kf\in X_0$ for each $k$, or ${\bf
R}$ linear and $~\mbox{}_kg(y)\in {\bf R}$ for each $k$ and every
$y\in U$, at each point $x\in U$.
\par When an operator $(I-A_x)$ is invertible and Condition $(RS)$
is satisfied, Equation $(12)$ can be resolved:
\par $(13)$ $\sum_k ~\mbox{}_kf(x,y)~\mbox{}_kg(y)=(I-A_x)^{-1}u(x,y)$.
\par If condition $(RS)$ is not fulfilled, the corresponding system
of equations in real components $(A_x)_{j,s}$, $~\mbox{}_kf_s$ and
$~\mbox{}_kg_s$ can be considered.
\par If $A: X\to X$ is a bounded ${\bf R}$ linear operator on a
Banach space $X$ with the norm $ \| A \| <1$, then \par $(I-A)^{-1}
= \sum_{n=0}^{\infty } A^n$. \\
The anti-derivative operator $g\mapsto ~\mbox{}_{\sigma }
\int_{\mbox{}_0x}^x g(z)dz$ is compact from $C^0(V,{\cal A}_r)$ into
$C^0(V,{\cal A}_r)$ for a compact domain $V$ in ${\cal A}_r$, where
$C^0(V,{\cal A}_r)$ is the Banach space over ${\cal A}_r$ of all
continuous functions $g: V\to {\cal A}_r$ supplied with the supremum
norm $\| g \| := \sup_{x\in V} |g(x)|$, $~\mbox{}_0x$ is a marked
point in $V$, $ ~ x\in V$. Therefore, an operator $A_x$ will be
invertible with a suitable choice of a function ${\sf F}$ satisfying
the system of ${\bf R}$ linear partial differential equations $(4)$
or $(5)$ and a real non-zero parameter ${\sf p}\ne 0$.

\par {\bf 7.1. Proposition.} {\it Let $V$ be a compact domain in the
Cayley-Dickson algebra ${\cal A}_r$ with $2\le r$  and let its
foliation be with rectifiable paths and satisfy Condition 2$(4)$ and
$\Lambda $ be a compact subset in ${\bf R}^{2^r}$. Then the
anti-derivative operator $\mbox{}_{\sigma }\int ^x$ from \S 2 is
compact from $C^0(V,{\cal A}_r)$ into $C^0(V,{\cal A}_r)$.}
\par {\bf Proof.} As usually $C^1(V,{\cal A}_r)$ denotes the space of all continuously
differentiable functions $f: V\to {\cal A}_r$ with the supremum norm
\par $\| f \| _{C^1} := \sup_{x\in V} |f(x)| + \sum_{j=0}^{2^r-1}
\sup_{x\in V} |\partial f(x)/\partial x_j|$. \par The decomposition
is valid: $C^s(V,{\cal A}_r)= C^s(V,{\bf R})i_0\oplus...\oplus
C^s(V,{\bf R})i_{2^r-1}$ for $s=0,1,...$. On the other hand, the
partial differential operators $C^1(V,{\bf R})\ni f\mapsto i_{\xi
(j)}^*(\partial f(x)/\partial x_j)\in C^0(V,{\bf R})$ are ${\bf R}$
linearly independent, since $Re (i_ji_k^*)=0$ for each $j\ne k$.
 In view of Theorems 2.6 \cite{ludfov,ludoyst} and 1.2.7, 4.2.5
and 4.2.23 and Corollary 4.2.6 \cite{ludancdnb} the anti-derivative
mapping $C^0(V,{\cal A}_r)\ni g\mapsto \mbox{}_{\sigma }\int
^xg(z)dz\in C^1(V,{\cal A}_r)$ is continuous, since $\sigma _x~
\mbox{}_{\sigma }\int^x g(z)dz=g(x)$. But the embedding $C^1(V,{\cal
A}_r)\hookrightarrow C^0(V,{\cal A}_r)$ is the ${\bf R}$ linear
${\cal A}_r$ additive compact operator. Therefore, the
anti-derivative operator $\mbox{}_{\sigma }\int ^x: C^0(V,{\cal
A}_r)\to C^0(V,{\cal A}_r)$ is compact.

\par {\bf 8. Example.} Let \par $(1)$ $N(x,z,y)={\sf K}(x,z)$ in 7$(1)$. \\ Acting from the
left on both sides of Equation 7$(1)$ one gets:
\par $(2)$ $L_s{\sf K}(x,y)=L_s {\sf p} ~ \mbox{}_{\sigma }
\int_x^{\infty } {\sf F}(z,y){\sf K}(x,z)dz $. \par Take the
hyperbolic partial differential operator \par $(3)$
$L_1=\mbox{}_1\sigma _x^2- ~ \mbox{}_2\sigma _y^2$ \\ of the second
order, where the operator $\mbox{}_k\sigma $ is with coefficients
$\mbox{}_k\psi _j\in {\bf R}$ for each $j=0,...,2^r-1$ and a
transposition $\xi _k\in S_{2^r}$ (see \S 1.4). Condition 7$(4)$ for
$s=1$ means that $L_1F=0$, i.e. \par $(4)$ $\mbox{}_1\sigma
_z^2F(z,y)=\mbox{}_2\sigma _y^2F(z,y)$. \\ Then due to Proposition 4
and Corollary 5 we deduce that
\par $(\mbox{}_1\sigma _x^2- ~ \mbox{}_2\sigma _y^2){\sf K}(x,y) =
{\sf p}(\mbox{}_1\sigma _x^2- ~ \mbox{}_2\sigma _y^2)  ~
\mbox{}_{\mbox{}_1\sigma } \int_x^{\infty } {\sf F}(z,y){\sf
K}(x,z)dz$
\par $ ={\sf p} ~ \mbox{}_1^2\sigma _x^2 ~\mbox{}_{\mbox{}_1\sigma }
\int_x^{\infty } {\sf F}(z,y){\sf K}(x,z)dz +{\sf p}~
\mbox{}_1A_2({\sf F},{\sf K})(x,y) - {\sf p}~ \mbox{}_1^1\sigma_z^2~
\mbox{}_{\mbox{}_1\sigma }\int_x^{\infty } {\sf F}(z,y){\sf
K}(x,z)dz$\par $={\sf p} (\mbox{}_1^2\sigma _x^2 - ~
\mbox{}_1^2\sigma _z^2) ~\mbox{}_{\mbox{}_1\sigma }\int_x^{\infty }
{\sf F}(z,y){\sf K}(x,z)dz+ {\sf p} ~ \mbox{}_1A_2({\sf F},{\sf
K})(x,y) - {\sf p} ~ \mbox{}_1B_2({\sf F},{\sf K})(x,y)$,
consequently,
\par $(5)$ $(I-{\sf A}_x)[(\mbox{}_1\sigma _x^2- ~ \mbox{}_2\sigma
_y^2){\sf K}(x,y)]={\sf p} ~ \mbox{}_1A_2({\sf F},{\sf K})(x,y)-
{\sf p}~ \mbox{}_1B_2({\sf F},{\sf K})(x,y)$\par $ = - 2{\sf
p}~\mbox{}^2_1\sigma _x[{\sf F}(x,y) {\sf K}(x,x)] = - 2{\sf
p}~\mbox{}^2_1\sigma _x\{ [(I-{\sf A}_x){\sf K}(x,y)]{\sf K}(x,x) \}
$, \\ where the terms $\mbox{}_kA_s$ and $\mbox{}_kB_s$ correspond
to $\mbox{}_k\sigma $.
\\ If the Cayley-Dickson algebra ${\cal A}_r$ is either with $r=
2$ and ${\sf F}\in Mat_s({\bf R})$ and ${\sf K}\in Mat_s({\cal
A}_2)$, $~ s\in {\bf N}$, or ${\sf F} \in {\bf R}$ and ${\sf K}\in
{\cal A}_3$ with $r=3$, then
\par $~\mbox{}^2_1 \sigma_x \{ [(I-{\sf A}_x){\sf K}(x,y)]{\sf K}(x,x) \} =
[(I-{\sf A}_x){\sf K}(x,y)] [~\mbox{}_1\sigma _x {\sf K}(x,x)] $, \\
since ${\bf R}$ is the center of the Cayley-Dickson algebra ${\cal
A}_r$ with $r\ge 2$, $ ~ <{\sf F}(z,y),{\sf K}(x,z),~\mbox{}_1\sigma
_x {\sf K}(x,x)>=0$ in these cases, where $<a,b,c> := (ab)c-a(bc)$
denotes the associator of Cayley-Dickson numbers $a, b, c\in {\cal
A}_r$. Then Conditions 7$(1,2)$ imply that the function ${\sf K}$
satisfies the non-linear partial differential equation:
\par $(6)$ $(\mbox{}_1\sigma _x^2- ~ \mbox{}_2\sigma _y^2){\sf
K}(x,y)+2{\sf p} {\sf K}(x,y)[~\mbox{}_1\sigma _x{\sf K}(x,x)]=0$.
\par If put $u(x)=2~\mbox{}_1\sigma _x {\sf K}(x,x)$ over the
quaternion skew field ${\bf H}={\cal A}_2$, i.e. for $r=2$, and
substitute ${\sf K}(x,y)=\Phi (x,k)\exp (J Re (ky))$ into $(6)$, we
deduce that a function $\Phi $ satisfies Schr\"odinger's equation:
\par $(7)$ $~\mbox{}_1\sigma _x^2\Phi (x,k) + \Phi (x,k)({\sf p}u +
\sum_jk_j^2i_j^2~\mbox{}_2\psi _j^2)=0$, \\
where $k\in {\bf H}$, since ${\bf H}$ is associative and the
generator $J$ commutes with $i_0,...,i_{2^r-1}$, also $\exp (J Re
(ky))$ commutes with $\Phi (x,k)$.

\par Now we take the third order partial differential operator with
$\mbox{}_1\psi _0= ~\mbox{}_2\psi _0 =0$:
\par $(8)$ $L_2f=(\mbox{}_3\sigma _t + ~ \mbox{}_1\sigma ^3_x+
~ 3 ~ \mbox{}_2\sigma _y ~ \mbox{}_1\sigma ^2_x + ~ 3 ~
\mbox{}_2\sigma ^2_y ~ \mbox{}_1\sigma _x + ~ \mbox{}_2\sigma ^3_y)
f$. Then we put
\par $(9)$ $L_{2,j}f=[\mbox{}_3\psi _{\xi _3(j)} \partial _{t_j} + ~
(\mbox{}_1\psi _j ~ \partial _{x_{\xi _1(j)}} + ~3~ \mbox{}_2\psi _j
\partial _{y_{\xi _2(j)}}) ~\mbox{}_1\sigma _x^2+ ( \mbox{}_2\psi
_j~ \partial _{y_{\xi _2(j)}}+ ~3~ \mbox{}_1\psi _j
\partial _{x_{\xi _1(j)}}) ~\mbox{}_2\sigma _y^2] f$
\\ and impose Condition 7$(5)$ for $s=2$:
\par $(10)$ $L_{2,j} {\sf F}(x,y) =0$ for each $j$, \\
where  the operator $~\mbox{}_3\sigma _t$ is with real constant
coefficients $\mbox{}_3\psi _j$ and a transposition $\xi _3\in
S_{2^r}$. We suppose that functions $\sf F$ and $\sf K$ may depend
on $t$. Therefore, we get from Equation $(2)$ and Condition $(10)$:
\par $(11)$ $(\mbox{}_3\sigma _t + ~ \mbox{}_1\sigma ^3_x+
~ 3 ~ \mbox{}_2\sigma _y ~ \mbox{}_1\sigma ^2_x + ~ 3 ~
\mbox{}_2\sigma ^2_y ~ \mbox{}_1\sigma _x + ~ \mbox{}_2\sigma
^3_y){\sf K}(x,y) $\par $= {\sf p} (\mbox{}_3\sigma _t + ~
\mbox{}_1\sigma ^3_x+ ~ 3  ~ \mbox{}_2\sigma _y ~ \mbox{}_1\sigma
^2_x + ~ 3  ~ \mbox{}_2\sigma ^2_y ~ \mbox{}_1\sigma _x + ~
\mbox{}_2\sigma ^3_y) ~ \mbox{}_{\mbox{}_1\sigma } \int_x^{\infty }
{\sf F}(z,y){\sf K}(x,z)dz=$
\par $ - {\sf p} (\mbox{}_1^1\sigma ^3_z+
~ 3 ~ \mbox{}_2^1\sigma _y ~ \mbox{}_1^1\sigma ^2_z + ~ 3 ~
\mbox{}_2^1\sigma ^2_y ~ \mbox{}_1^1\sigma _z + ~ \mbox{}_2^1\sigma
^3_y)~ \mbox{}_{\mbox{}_1\sigma }\int_x^{\infty } {\sf F}(z,y){\sf
K}(x,z)dz +$
\par ${\sf p} (\mbox{}_3^2\sigma _t + ~\mbox{}_1\sigma ^3_x+ ~ 3 ~ \mbox{}_2\sigma _y ~
\mbox{}_1\sigma ^2_x + ~ 3 ~ \mbox{}_2\sigma ^2_y ~ \mbox{}_1\sigma
_x + ~ \mbox{}_2\sigma ^3_y) ~ \mbox{}_{\mbox{}_1\sigma }
\int_x^{\infty } {\sf F}(z,y){\sf K}(x,z)dz=I+ {\sf p} ~
\mbox{}_3^2\sigma _t~
\mbox{}_{\mbox{}_1\sigma } \int_x^{\infty } {\sf F}(z,y){\sf K}(x,z)dz$, \\
where $I=I_1+I_2+I_3$,
\par $(12)$ $I_1= {\sf p} (\mbox{}_1\sigma ^3_x
- ~ \mbox{}_1^1\sigma ^3_z) ~ \mbox{}_{\sigma } \int_x^{\infty }
{\sf F}(z,y){\sf K}(x,z)dz$
\par $={\sf p} (\mbox{}_1^2\sigma ^3_x
+ ~ \mbox{}_1^2\sigma ^3_z) ~ \mbox{}_{\mbox{}_1\sigma }
\int_x^{\infty } {\sf F}(z,y){\sf K}(x,z)dz + {\sf p}~
\mbox{}_1A_3({\sf F},{\sf K})(x,y) - {\sf p}~ \mbox{}_1B_3({\sf
F},{\sf K})(x,y) $
\par $(13)$ $I_2 = 3{\sf p} ( \mbox{}_2\sigma _y ~ \mbox{}_1\sigma ^2_x - ~
\mbox{}_2\sigma _y ~ \mbox{}_1^1\sigma ^2_z ) ~
\mbox{}_{\mbox{}_1\sigma } \int_x^{\infty } {\sf F}(z,y){\sf
K}(x,z)dz$
\par $=~ 3 {\sf p}~ \mbox{}_2\sigma _y [(\mbox{}_1^2\sigma ^2_x - ~
\mbox{}_1^2\sigma ^2_z )~ \mbox{}_{\mbox{}_1\sigma } \int_x^{\infty
} {\sf F}(z,y){\sf K}(x,z)dz + ~\mbox{}_1A_2({\sf F},{\sf K})(x,y) -
~\mbox{}_1B_2({\sf F},{\sf K})(x,y)]$,
\par $(14)$ $I_3 = 3{\sf p}( ~ \mbox{}_2\sigma ^2_y ~ \mbox{}_1\sigma _x - ~
\mbox{}_2\sigma ^2_y ~ \mbox{}_1^1\sigma _z ) ~
\mbox{}_{\mbox{}_1\sigma } \int_x^{\infty } {\sf F}(z,y){\sf
K}(x,z)dz$
\par $= 3 {\sf p}~ \mbox{}_2\sigma ^2_y (~ \mbox{}_1^2\sigma _x + ~
\mbox{}_1^2\sigma _z ) ~ \mbox{}_{\mbox{}_1\sigma } \int_x^{\infty }
{\sf F}(z,y){\sf K}(x,z)dz$, \\ since $\mbox{}_1A_1=\mbox{}_1B_1$.
Equations $(4,14)$, 4$(3)$ and 5$(1,3)$ imply that
\par $(15)$ $I_3= 3 {\sf p}~ \mbox{}_1^1\sigma ^2_z (~ \mbox{}_1^2\sigma _x + ~
\mbox{}_1^2\sigma _z ) ~ \mbox{}_{\mbox{}_1\sigma } \int_x^{\infty }
{\sf F}(z,y){\sf K}(x,z)dz$ \\ $=3 {\sf p}(\mbox{}_1^2\sigma _z - ~
\mbox{}_1^1\sigma _x) \{ (~ \mbox{}_1^2\sigma _x + ~
\mbox{}_1^2\sigma _z )[{\sf F}(x,y){\sf K}(x,z)] \} |_{z=x} + ~ 3
{\sf p}(\mbox{}_1^2\sigma _x ~ \mbox{}_1^2\sigma _z^2 + ~
\mbox{}_1^2\sigma ^3_z )~ \mbox{}_{\mbox{}_1\sigma } \int_x^{\infty
} {\sf F}(z,y){\sf K}(x,z)dz.$ \\ Particularly, if the
Cayley-Dickson algebra ${\cal A}_r$ is with $2\le r\le 3$ and ${\sf
F} \in Mat_s({\bf R})$ and ${\sf K}\in Mat_s({\cal A}_r)$, $ ~s\in
{\bf N}$ for $r=2$, $~ s=1$ for $r=3$, then Equations $(6,13)$ imply
that
\par $(16)$ $I_2= -3 ~ \mbox{}_2\sigma _y [{\sf K}(x,y) - {\sf F}(x,y)]
u(x) + ~ 3{\sf p} ~ \mbox{}_2\sigma _y [\mbox{}_1A_2({\sf F},{\sf
K})(x,y) - ~\mbox{}_1B_2({\sf F},{\sf K})(x,y)]$, \\
since
\par $(17)$ $ -2 {\sf p} ~\mbox{}_{\mbox{}_1\sigma } \int _x^{\infty } {\sf
F}(z,y) [{\sf K}(x,y) (\mbox{}_1\sigma _x {\sf K}(x,x))]dz = - {\sf
p} ~\mbox{}_{\mbox{}_1\sigma } \int _x^{\infty } {\sf F}(z,y) [{\sf
K}(x,z) u(x)]dz$
\par $=- {\sf p} [~\mbox{}_{\mbox{}_1\sigma } \int _x^{\infty } {\sf F}(z,y)
{\sf K}(x,z)dz] u(x) =[{\sf K}(x,y)-{\sf F}(x,y)]u(x)$, \\
 where $u(x) =2~\mbox{}_1\sigma _x {\sf K}(x,x)$. We deduce
from Formulas $(6,10,11,15-17)$ that
\par $(18)\quad (\mbox{}_3\sigma _t + ~ \mbox{}_1\sigma ^3_x+
~ 3 ~ \mbox{}_2\sigma _y ~ \mbox{}_1\sigma ^2_x + ~ 3 ~
\mbox{}_2\sigma ^2_y ~ \mbox{}_1\sigma _x + ~ \mbox{}_2\sigma
^3_y){\sf K}(x,y) + 3\mbox{}_2^1\sigma _y[{\sf K}(x,y)u(x)]=$
\par ${\sf p} (\mbox{}_3^2\sigma _t + ~ \mbox{}_1^2\sigma ^3_x+ ~ 3 ~
\mbox{}_1^2\sigma _z ~ \mbox{}_1^2\sigma ^2_x + ~ 3 ~
\mbox{}_1^2\sigma ^2_z ~ \mbox{}_1^2\sigma _x + ~ \mbox{}_1^2\sigma
^3_z)~ \mbox{}_{\mbox{}_1\sigma } \int_x^{\infty } {\sf F}(z,y){\sf
K}(x,z)dz$\par $ +3{\sf p}~ \mbox{}_{\mbox{}_1\sigma }
\int_x^{\infty } {\sf F}(z,y)[\mbox{}^1_2\sigma _z {\sf
K}(x,z)u(x)]dz  +T$, where \par $T={\sf p} ~\mbox{}_1A_3({\sf
F},{\sf K})(x,y) - {\sf p} ~\mbox{}_1B_3({\sf F},{\sf K})(x,y) + ~3
{\sf p} ~ \mbox{}_2 \sigma _y [\mbox{}_1A_2({\sf F},{\sf K})(x,y) -
~ \mbox{}_1B_2({\sf F},{\sf K})(x,y)] +$\par $ ~ 3 {\sf p} ~
\mbox{}_2^1\sigma _y[F(x,y)u(x)] + 3{\sf p} ~(~\mbox{}_1^2\sigma _z
~ \mbox{}_1^2\sigma _x + ~ \mbox{}_1^2\sigma _z^2 - ~
\mbox{}_1^1\sigma _x ~ \mbox{}_1^2\sigma _x - \mbox{}_1^1\sigma _x ~
\mbox{}_1^2\sigma _z) [{\sf F}(x,y){\sf K}(x,z)]|_{z=x}$.
\par Then for each two continuously differentiable ${\cal A}_r$ valued functions
${\sf G}(z,y)$ and ${\sf K}(x,z)$ one has \par $(\mbox{}^1\sigma _z-
~ \mbox{}^2\sigma _x - ~\mbox{}^2\sigma _z)[{\sf G}(z,y){\sf
K}(x,z)]= (\mbox{}^1\sigma _z+ ~\mbox{}^2\sigma _x +
~\mbox{}^2\sigma _z) [{\sf G}(z,y)\check{{\sf K}}(x,z)]= (\sigma _x+
\sigma _z) [{\sf G}(z,y)\check{{\sf K}}(x,z)]$, \\ where
$\check{{\sf K}}(x,z) := {\sf K}(-x,-z)$ for each $-x, -z\in U$.
Therefore, the identity
\par $[\mbox{}_{\sigma }\int_x ,[ ~\mbox{}^1\sigma _x,
~\mbox{}^2\sigma _x]] \{ {\sf G}(z,y){\sf K}(x,z) \} dz = 0 $ \\ is
satisfied, consequently, in the considered case $2\le r\le 3$ and
${\sf F} \in Mat_s({\bf R})$ and ${\sf K}\in Mat_s({\cal A}_r)$, $
~s\in {\bf N}$ for $r=2$, $~ s=1$ for $r=3$, we get
\par $(19)\quad [(I-A_x),[\mbox{}^1\sigma _x,\mbox{}^2\sigma _x]] \{ {\sf
F}(x,y){\sf K}(x,x) \} =0$, \\ since ${\bf R}$ is the center of the
Cayley-Dickson algebra ${\cal A}_r$ and \par $(20)\quad
[\mbox{}^2\sigma _x, ~ \mbox{}^3\sigma _x] {\sf F}(z,y)\{ {\sf
F}(x,z){\sf K}(x,x) \} = {\sf F}(z,y) ([\mbox{}^1\sigma _x, ~
\mbox{}^2\sigma _x] \{ {\sf F}(x,z){\sf K}(x,x) \} )$ \\ due to
Formulas $(1)$, 4$(7,8)$ and 7$(2)$, since $[\mbox{}^1\sigma _x+
~\mbox{}^2\sigma _x, \mbox{}^1\sigma _x - ~\mbox{}^2\sigma _x] = - 2
[\mbox{}^1\sigma _x ,  ~\mbox{}^2\sigma _x]\{ {\sf G}(x,y){\sf
K}(x,x) \} $.

\par Substituting the expressions of $A_2-B_2$ and $A_3-B_3$ from
Corollary 5 and using Formulas $(2,6,19,20)$ one gets:
\par $(21)\quad T=- {\sf p} (3 ~\mbox{}^2_1\sigma ^2_x + ~ \mbox{}^2_1\sigma _x ~ \mbox{}^2_1\sigma
_z+ 2~ \mbox{}^2_1\sigma _z ~ \mbox{}^2_1\sigma _x) [{\sf
F}(x,y){\sf K}(x,z)]|_{z=x}$\par $ - {\sf p} (2 ~\mbox{}^1_1\sigma
_x ~\mbox{}^2_1\sigma _x + ~\mbox{}^2_1\sigma _x ~\mbox{}^1_1\sigma
_x) [{\sf F}(x,y){\sf K}(x,x)]$\par $+ 3(1-{\sf p}) ~\mbox{}_2\sigma
_y [{\sf F}(x,y)u(x)] + 3{\sf p} ~(~\mbox{}_1^2\sigma _z ~
\mbox{}_1^2\sigma _x + ~ \mbox{}_1^2\sigma _z^2 - ~
\mbox{}_1^1\sigma _x ~ \mbox{}_1^2\sigma _x - \mbox{}_1^1\sigma _x ~
\mbox{}_1^2\sigma _z) [{\sf F}(x,y){\sf K}(x,z)]|_{z=x} $

\par $= - 3{\sf p} ~\mbox{}_1^1\sigma _x[{\sf F}(x,y)u(x)] -
3{\sf p} (\mbox{}_1^2\sigma _x^2- ~ \mbox{}_1^2\sigma _z^2)[{\sf
F}(x,y){\sf K}(x,z)]|_{z=x} + 3(1-{\sf p}) ~\mbox{}_2\sigma _y [{\sf
F}(x,y)u(x)] $

\par $+ {\sf p} [ \mbox{}^2_1\sigma _z, ~ \mbox{}^2_1\sigma _x]
[{\sf F}(x,y) {\sf K}(x,z)]|_{z=x} + {\sf p} [ \mbox{}^1_1\sigma _x,
~ \mbox{}^2_1\sigma _x] [{\sf F}(x,y) {\sf K}(x,x)]$
\par $= - 3{\sf p} ~\mbox{}_1^1\sigma _x[{\sf F}(x,y)u(x)] + 3{\sf p} {\sf F}(x,y)[{\sf
K}(x,x)u(x)]+ 3(1-{\sf p}) ~\mbox{}_2\sigma _y [{\sf F}(x,y)u(x)] $
\par $+ {\sf p} [ \mbox{}^2_1\sigma _z, ~ \mbox{}^2_1\sigma _x]
[{\sf F}(x,y) {\sf K}(x,z)]|_{z=x} + {\sf p} [ \mbox{}^1_1\sigma _x,
~ \mbox{}^2_1\sigma _x] [{\sf F}(x,y) {\sf K}(x,x)]$

\par $=- 3{\sf p} ~\mbox{}_1^1\sigma _x[{\sf K}(x,y)u(x)] + 3{\sf p}^2 ~\mbox{}_1\sigma
_x\{ (\mbox{}_{\mbox{}_1\sigma }\int_x^{\infty } {\sf F}(z,y){\sf
K}(x,z)dz)u(\eta ) \}|_{\eta =x} + 3{\sf p} {\sf F}(x,y)[{\sf
K}(x,x)u(x)]$ \par $ + 3(1-{\sf p}) ~\mbox{}_2\sigma _y [{\sf
F}(x,y)u(x)] + {\sf p} [ \mbox{}^2_1\sigma _z, ~ \mbox{}^2_1\sigma
_x] [{\sf F}(x,y) {\sf K}(x,z)]|_{z=x} + {\sf p} [ \mbox{}^1_1\sigma
_x, ~ \mbox{}^2_1\sigma _x] [{\sf F}(x,y) {\sf K}(x,x)]$
\par $= - 3{\sf p} ~\mbox{}_1^1\sigma _x[{\sf F}(x,y)u(x)] + 3{\sf p}^2 ~\mbox{}_1^2\sigma
_x[(\mbox{}_{\mbox{}_1\sigma }\int_x^{\infty } [({\sf F}(z,y){\sf
K}(x,z))u(x)]dz $
\par $ + 3(1-{\sf p}) ~\mbox{}_2\sigma _y [(I-A_x){\sf
K}(x,y)u(x)] + {\sf p}  [(I-A_x)  {\sf K}(x,y)] \{ [ \mbox{}_1\sigma
_z, ~ \mbox{}_1\sigma _x] {\sf K}(x,z)]|_{z=x} \} + {\sf p}(I-A_x)
\{ [ \mbox{}^1_1\sigma _x, ~ \mbox{}^2_1\sigma _x] [{\sf K}(x,y)
{\sf K}(x,x)] \} $, \\
since $2\le r\le 3$, $F \in Mat_s({\bf R})$ and ${\sf K}\in
Mat_s({\cal A}_r)$ with $s\in {\bf N}$ for $r=2$ and $s=1$ for
$r=3$, hence $<{\sf F}(x,y), {\sf K}(x,x), u(x)>=0$, where $<a,b,c>
:= (ab)c-a(bc)$ denotes the associator for each Cayley-Dickson
numbers $a, b, c \in {\cal A}_r$, since $i_0i_k-i_ki_0=0$,
$~i_ji_k=-i_ki_j$ for each $j\ne k\ge 1$,
\par $(22)$ $[\sigma _z,\sigma _x] = (\sigma
_z\sigma _x-\sigma _x\sigma _z) = - [\sigma _x,\sigma _z]$.\\
We now take ${\sf p} =1$. Therefore, in accordance with Formulas
$(19-21)$ and 7$(13)$ the equality
\par $(23)$ $(\mbox{}_3\sigma _t + ~ \mbox{}_1\sigma ^3_x+ ~ 3  ~
\mbox{}_2\sigma _y ~ \mbox{}_1\sigma ^2_x + ~ 3 ~ \mbox{}_2\sigma
^2_y ~ \mbox{}_1\sigma _x + ~ \mbox{}_2\sigma ^3_y){\sf K}(x,y)
$\par $ + 6 (\mbox{}_1^1\sigma _x + ~ \mbox{}_2^1\sigma _y) [{\sf
K}(x,y)(~\mbox{}_1\sigma _x {\sf K}(x,x))] - {\sf K}(x,y) \{ [
\mbox{}_1\sigma _z, ~ \mbox{}_1\sigma _x] {\sf K}(x,z)]|_{z=x}
\}$\par $ - [ \mbox{}^1_1\sigma _x, ~ \mbox{}^2_1\sigma _x] [{\sf
K}(x,y) {\sf K}(x,x)] =0$ \\ follows, when the operator $(I-A_x)$ is
invertible. Particularly, for $s=1$ and $\mbox{}_1\sigma
=\mbox{}_2\sigma $ with $\mbox{}_1\psi _0=0$ and $\mbox{}_3\sigma
_t=\partial /\partial t_0$ differentiating Equation $(23)$ with the
operator $\mbox{}_1\sigma _x$ and then restricting on the diagonal
$x=y$ and taking into account Formulas $(1,6)$, 7$(1)$ and $(19)$
one gets the equation
\par $(24)$ $u_t(t,x)+6~\mbox{}_1\sigma _x[u(t,x)u(t,x)] +
~\mbox{}_1\sigma _x^3u(t,x)=0$\\
of Korteweg-de-Vries' type, where $u(t,x)=2~\mbox{}_1\sigma _x{\sf
K}(x,x)$, since $\sigma (f(x)g(x))=(\mbox{}^1\sigma + ~
\mbox{}^2\sigma )(f(x)g(x))$, $~[\mbox{}^1\sigma , ~ \mbox{}^2\sigma
](f(x)g(x))= - \frac{1}{2} [\sigma , ~\mbox{}^1\sigma - ~
\mbox{}^2\sigma ](f(x)g(x))$, $~\sigma [\mbox{}^1\sigma , ~
\mbox{}^2\sigma ](f(x)g(x))= - [\mbox{}^1\sigma , ~ \mbox{}^2\sigma
]\sigma (f(x)g(x))$ and hence $~\sigma ^{-1} [\mbox{}^1\sigma , ~
\mbox{}^2\sigma ](f(x)g(x)) = - [\mbox{}^1\sigma , ~ \mbox{}^2\sigma
]\sigma ^{-1} (f(x)g(x))$ for $\sigma =\mbox{}_1\sigma _x$ with
$\mbox{}_1\psi _0=0$ and continuously differentiable functions
$f(x)$ and $g(x)$ (see Remark 3.2).

\par {\bf 8.1. Theorem.} {\it A solution of partial differential Equation
$(23)$ with $\mbox{}_1\psi _0=\mbox{}_2\psi _0=0$ over the
Cayley-Dickson algebra ${\cal A}_r$ with $2\le r\le 3$ is given by
Formulas $(2-4,9,10)$ with ${\sf p} =1$ whenever the appearing
integrals uniformly converge by parameters on compact sub-domains
(see Proposition 4) and the operator $(I-A_x)$ is invertible and
${\sf F}\in Mat_s({\bf R})$ and ${\sf K}\in Mat_s({\cal A}_r)$,
$s\in {\bf N}$ when $r=2$, $s=1$ when $r=3$.}

\par {\bf 9. Example.} Consider the integral equation
\par $(1)$ ${\sf K}(x,y) = {\sf F}(x,y) + \frac{\sf p}{4} ~
\mbox{}_{\sigma } \int_x^{\infty }(~ \mbox{}_{\sigma
}\int_x^{\infty } {\sf F}(u,y)[{\sf F}(z,u){\sf K}(x,z)]dz)du$, \\
where ${\sf p}$ is a non-zero real constant. We take, for example,
\par $(2)$ $L_1= \sigma _x - ~\mbox{}_1\sigma _y$ with $\psi _0= ~\mbox{}_1\psi _0=0$.
\par A solution of the equation \par $(3)$ $L_1{\sf F}(x,y)=0$ \\ has the form
${\sf F}(x,y)={\sf G}(\frac{(a,x>+(\mbox{}_1a,y>}{2})$ or we shall
write ${\sf F}(\frac{(a,x>+(\mbox{}_1a,y>}{2})$ instead of ${\sf
F}(x,y)$, where \par $(4)$ $(a,x> := \sum_{j=0}^{2^r-1} a_jx_ji_j$ \\
for $a, x\in {\cal A}_r$, $ ~ x=\sum_{j=0}^{2^r-1} x_ji_j$, $ ~
x_j\in {\bf R}$ for each $j$, $~a_j= \psi _j$ for $\psi _j\ne 0$ and
$a_j=1$ for $\psi _j=0$. Using suitable change of real variables
$x_j, y_j$ we can suppose without loss of generality, that $\psi _j,
~\mbox{}_1\psi _j\in \{ 0, 1 \} $ for each $j=0,...,2^r-1$.
Therefore, a solution of Equation $(3)$ can be written as \par
$(3.1)$ ${\sf F}(\frac{(b,x+y>}{2})$ as well, where
$b_j=a_j~\mbox{}_1a_j\in \{ 0,1 \} $ for each $j$. Differentiable
functions of the form ${\sf F}(\frac{(b,x+y>}{2})$ are also
solutions of the system of partial differential equations
\par $(3.2)$ $L_{1,j}{\sf F}(\frac{(b,x+y>}{2})=0,$ where
\par $(3.3)$ $L_{1,j} = \psi _j \partial _{x_j}- ~\mbox{}_1\psi _j\partial _{y_j}.$

\par It is convenient to introduce the notation
\par $(5)$ ${\sf K}_2(x,z) := \mbox{}_{\sigma
}\int_0^{\infty } {\sf F}(\frac{(a,x+\zeta > +
(\mbox{}_1a,z>}{2}){\sf K}(x,x+\zeta )d\zeta $. \\
Using Condition $(3)$ we rewrite Equation $(1)$ in the form:
\par $(6)$ ${\sf K}(x,y) = {\sf
F}(\frac{(a,x>+(\mbox{}_1a,y>}{2})$\par $ + \frac{\sf p}{4}
~\mbox{}_{\sigma }\int_0^{\infty } ~\mbox{}_{\sigma }\int_0^{\infty
} [{\sf F}(\frac{(a,x+\eta > + (\mbox{}_1a,y>}{2}) ({\sf
F}(\frac{(a,x+\zeta > + (\mbox{}_1a,x+\eta >}{2}){\sf K}(x,x+\zeta
))d\zeta ]d\eta $\par $= {\sf F}(\frac{(a,x>+(\mbox{}_1a,y>}{2}) +
\frac{\sf p}{4} ~\mbox{}_{\sigma }\int_0^{\infty }{\sf
F}(\frac{(a,x+\eta > + (\mbox{}_1a,y>}{2}) {\sf K}_2(x,x+\eta )d\eta
$. \\ Then let us put:
\par $(7)$ $A_xf(y) := \frac{\sf p}{4}~\mbox{}_{\sigma }\int_0^{\infty } ~\mbox{}_{\sigma }\int_0^{\infty
} [{\sf F}(\frac{(a,x+\eta > + (\mbox{}_1a,y>}{2}) ({\sf
F}(\frac{(a,x+\zeta > + (\mbox{}_1a,x+\eta >}{2})f(x+\zeta )) d\zeta
]d\eta $ \\ for a continuous function $f$, consequently,
\par $(8)$ $(I-A_x){\sf K}(x,y) ={\sf
F}(\frac{(a,x>+(\mbox{}_1a,y>}{2})$. \par Particularly, if ${\sf F}
\in Mat_s({\bf R})$ and ${\sf K}\in Mat_s({\cal A}_r)$, then
\par $(9)$ $(I-A_x){\sf K}_2(x,y) = ~\mbox{}_{\sigma }\int_0^{\infty } {\sf
F}(\frac{(a,x+\zeta >+(\mbox{}_1a,y>}{2}) {\sf
F}(\frac{(a,x>+(\mbox{}_1a,x+\zeta >}{2}) d\zeta $. \\ Acting on
both sides of Equation $(1)$ with the operator $L_1$ leads to the
relation:
\par $(10)$ $(\sigma _x - ~\mbox{}_1\sigma _y) {\sf K}(x,y) =
\frac{\sf p}{4} ~ ~\mbox{}^2\sigma _x ~\mbox{}_{\sigma
}\int_0^{\infty } {\sf F}(\frac{(a,x+\eta >+(\mbox{}_1a,y>}{2}) {\sf
K}_2(x,x+\eta )d\eta $. \\ Now we deduce an expression for $\sigma
_x{\sf K}_2(x,x+\eta )$. From the definition of ${\sf K}_2$ the
identity
\par $(11)$ $(\sigma _x+~\mbox{}_1\sigma _z){\sf K}_2(x,z) = (\mbox{}^1\sigma _x+ ~\mbox{}^1_1\sigma _z)
~\mbox{}_{\sigma }\int_0^{\infty } {\sf F}(\frac{(a,x+\zeta
>+(\mbox{}_1a,z>}{2}) {\sf K}(x,x+\zeta )d\zeta $ \par $+
\mbox{}^2\sigma _x ~\mbox{}_{\sigma }\int_0^{\infty } {\sf
F}(\frac{(a,x+\zeta >+(\mbox{}_1a,z>}{2}) {\sf K}(x,x+\zeta )d\zeta
$ \\
follows. Using $(3.1)$ we get that \par $(12)$ $(\sigma _x+
~\mbox{}_1\sigma _z) {\sf F}(\frac{(a,x+\zeta >+(\mbox{}_1a,z>}{2})
=(\sigma _x+ ~\mbox{}_1\sigma _z) {\sf F}(\frac{(b,x+\zeta +z>}{2})
= \sigma _{\xi } {\sf F}(\xi )|_{\xi =\frac{(b,x+\zeta +z>}{2}}$, \\
since $b_ja_j=b_j$ and $b_j~\mbox{}_1a_j=b_j$ for each
$j=0,...,2^r-1$. Therefore, in accordance with Proposition 4 the
equality
\par $(13)$ $(\sigma _x+\sigma _z){\sf K}_2(x,z) = (\mbox{}^2\sigma _1- ~\mbox{}^2\sigma _2)
~\mbox{}_{\sigma }\int_0^{\infty } {\sf F}(\frac{(b,x+\zeta +z>}{2})
{\sf K}(x,x+\zeta )d\zeta -2{\sf F}(\frac{(b,x+z>}{2}){\sf K}(x,x)$
\\ follows, where $\sigma _1{\sf K}(x,z) := \sigma _x{\sf K}(x,z)$ and $\sigma
_2{\sf K}(x,z) := \sigma _z{\sf K}(x,z)$. We seek a solution ${\sf
K}(x,y)$ satisfying the condition:
\par $(14)$ $\sigma _y{\sf K}(x,y) = ~ \mbox{}_1\sigma _y{\sf
K}(x,y)$. \\ From Equations $(10,13,14)$ we get, that
\par $(15)$ $(\sigma _x+\sigma _z){\sf K}_2(x,z) = -2{\sf F}(\frac{(b,x+z>}{2}){\sf K}(x,x)$
\par $+\frac{\sf p}{4} ~ ~\mbox{}^3\sigma _x ~\mbox{}_{\sigma }\int_0^{\infty }
{\sf F}(\frac{(b,x+\zeta +z>}{2}) \{ ~\mbox{}_{\sigma
}\int_0^{\infty }[ {\sf F}(\frac{(b,2x+\eta +\zeta >}{2}) {\sf
K}_2(x,x+\eta )]d\eta \} d\zeta $, consequently, \par $(16)$
$(I-A_x)[(\sigma _x+\sigma _z){\sf K}_2(x,z)]=-2{\sf
F}(\frac{(b,x+z>}{2}){\sf K}(x,x)=-2 [(I-A_x){\sf K}(x,z)]{\sf
K}(x,x)$. \par If either $r=2$ or ${\sf F}\in Mat_s({\bf R})$ and
${\sf K}\in Mat_s({\cal A}_r)$, then one obtains from $(3.3,5,6,13)$
analogously the identities:
\par $(17)$ $(\sigma _x- ~\mbox{}_1\sigma _y) {\sf K}(x,y) =
\frac{\sf p}{4} ~ (\mbox{}^3\sigma _1 - ~\mbox{}^3\sigma _2)
~\mbox{}_{\sigma }\int_0^{\infty } {\sf F}(\frac{(b,x+\zeta +y>}{2})
\{ ~\mbox{}_{\sigma }\int_0^{\infty }[ {\sf F}(\frac{(a,2x+\eta
+\zeta >}{2}) {\sf K}(x,x+\zeta )]d\zeta - 2{\sf F}(\frac{(a,2x+\eta
>}{2}){\sf K}(x,x)\} d\eta $\par $=A_x[(\sigma _x-\sigma _y){\sf K}(x,y)]
- \frac{\sf p}{2} ~\mbox{}_{\sigma }\int_0^{\infty } {\sf
F}(\frac{(b,x+\zeta +y>}{2}) [ {\sf F}(\frac{(a,2x+\eta >}{2}) {\sf
K}(x,x)]d\eta $.\\  Therefore,
\par $(18)$ $(I-A_x) [(\sigma _x-\sigma _y){\sf K}(x,y)]= -
\frac{\sf p} {2} [(I-A_x){\sf K}_2(x,y)]{\sf K}(x,x)+ \frac{\sf p}
{2}~ \mbox{}_{\sigma }\int_0^{\infty } <{\sf F}(\frac{(b,x+\eta +y
>}{2}), ~{\sf F}(\frac{(b,2x+\eta >}{2}),~ {\sf K}(x,x ) > d\eta $.
\\ In the associative quaternion case ${\cal A}_2={\bf H}$ or when
${\sf F}\in Mat_s({\bf R})$ and ${\sf K}\in Mat_s({\cal A}_r)$ the
last additive in Formula $(18)$ vanishes, i.e. \par $(18.1)$
$(I-A_x) [(\sigma _x-\sigma _y){\sf K}(x,y)]= - \frac{\sf p} {2}
[(I-A_x){\sf K}_2(x,y)]{\sf K}(x,x)$. \\ On the other hand, one has
the identities:
\par $(19)$ ${\sf K}_2(x,x+\zeta ) = \mbox{}_{\sigma }\int_0^{\infty
} {\sf F}(\frac{(b,2x+\xi +\zeta >}{2}){\sf K}(x,x+\xi ) d\xi $ and

\par $(19.1)$ $(I-A_x){\sf K}_2(x,z) = \mbox{}_{\sigma }\int_0^{\infty
} {\sf F}(\frac{(b,x+\zeta +z>}{2}){\sf K}(x,x+\zeta ) d\zeta$\par
$-\frac{\sf p}{4} ~  \mbox{}_{\sigma }\int_0^{\infty } {\sf
F}(\frac{(b,x+z +\eta >}{2})\{ {\sf F}(\frac{(b,2x+\eta +\zeta
>}{2})[~  \mbox{}_{\sigma }\int_0^{\infty } {\sf
F}(\frac{(b,2x+\xi +\zeta >}{2}) {\sf K}(x,x+\xi ) d\xi ] d\zeta \}
d\eta =$\par $~  \mbox{}_{\sigma }\int_0^{\infty } {\sf
F}(\frac{(b,x+z +\zeta >}{2}) [(I-A_x) {\sf K}(x,x+\zeta )] d\zeta
$\par $ = ~  \mbox{}_{\sigma }\int_0^{\infty } {\sf F}(\frac{(b,x+z
+\zeta >}{2}) {\sf F}(\frac{(b,2x+\zeta >}{2}) d\zeta $, that is
\par $(20)$ $(I-A_x){\sf K}_2(x,z)= ~  \mbox{}_{\sigma }\int_0^{\infty } {\sf F}(\frac{(b,x+z
+\zeta >}{2}) {\sf F}(\frac{(b,2x+\zeta >}{2}) d\zeta $. \\
If either $r=2$, i.e. ${\cal A}_2={\bf H}$, or when ${\sf F}\in
Mat_s({\bf R})$ and ${\sf K}\in Mat_s({\cal A}_r)$, $~s\in {\bf N}$
for $r=2$, $~s=1$ for $r=3$, Equations $(16, 18.1)$ imply that the
functions ${\sf K}$ and ${\sf K}_2$ are solutions of the system of
partial differential equations:
\par $(21)$ $(\sigma _x+\sigma _z){\sf K}_2(x,z) = -2{\sf K}(x,z){\sf
K}(x,x)$ and
\par $(22)$ $(\sigma _x-\sigma _z){\sf K}(x,z) = -\frac{\sf p}{2} {\sf
K}_2(x,z){\sf K}(x,x)$. \\ The action on both sides of Equation
$(1)$ with the operator $\sigma _x+ ~\mbox{}_1\sigma _y$ leads to
the identities:
\par $(23)$ $(\sigma _x + ~ \mbox{}_1\sigma _y) {\sf K}(x,y) =
\sigma _zF(z)|_{z=(b,x+y>/2} $\par $+\frac{\sf p}{4}
 ~ \mbox{}^3\sigma _x ~ \mbox{}_{\sigma
}\int_0^{\infty } ~ \mbox{}_{\sigma }\int_0^{\infty } \{ {\sf
F}(\frac{(b,x+\eta +y
>}{2}) [{\sf F}(\frac{(b,2x+\zeta +\eta >}{2}) {\sf K}(x,x+\zeta
)]d\zeta \} d\eta $\par $- \frac{\sf p}{2} {\sf
F}(\frac{(b,x+y>}{2})~ \mbox{}_{\sigma }\int_0^{\infty } [{\sf
F}(\frac{(b,2x+\zeta
>}{2}) {\sf K}(x,x+\zeta )]d\zeta $, \\ since
\par $(\mbox{}^1\sigma _x+ ~\mbox{}_1^1\sigma _y+ ~\mbox{}^2\sigma _x+ ~\mbox{}_1^2\sigma _y)
\{ {\sf F}(\frac{(b,x+\eta +y >}{2}) [{\sf F}(\frac{(b,2x+\zeta
+\eta >}{2}) g(x,\zeta ) ] \} $\par $ = 2\sigma _{\eta } \{ {\sf
F}(\frac{(b,x+\eta +y >}{2}) [{\sf F}(\frac{(b,2x+\zeta +\eta >}{2})
g(x,\zeta ) ] \} $ \\ for each ${\cal A}_r$ valued function
$g(x,\zeta )$, consequently,
\par $(24)$ $(I-A_x) [(\sigma _x+ ~\mbox{}_1\sigma _y){\sf
K}(x,y)]=\sigma _zF(z)|_{z=(b,x+y>/2} - \frac{\sf p}{2} [(I-A_x){\sf
K}(x,y)]{\sf K}_2(x,x)$. \par Then we use the condition \par $(25)$
$L_{2,j}{\sf F}=0$ for each $j=0,...,2^r-1$,
\\ where
\par $(26)$ $L_{2,j} = ~\mbox{}_2\psi _j\partial _{t_j}+\psi _j
\partial _{x_j}\sigma _x^2 +3 ~\mbox{}_1\psi _j\partial_{y_j} \sigma _x^2 + 3
\psi _j ~\mbox{}_1\sigma _y^2 \partial _{x_j} + ~\mbox{}_1\psi
_j\partial_{y_j} ~\mbox{}_1\sigma _y^2$ \\
and act on both sides of Equation $(1)$ with the operator
\par $L_2=\mbox{}_2\sigma _t +\sigma _x^3 +3 ~\mbox{}_1\sigma _y\sigma
_x^2 + 3 ~\mbox{}_1\sigma _y^2 \sigma _x + ~\mbox{}_1\sigma _y^3$
with $\psi _0=0$, $~\mbox{}_1\psi _0=0$, that gives
\par $(27)$ $L_2{\sf K}(x,y) =\frac{\sf p}{4}
(\mbox{}_2\sigma _t +\sigma _x^3 +3 ~\mbox{}_1\sigma _y\sigma _x^2 +
3 ~\mbox{}_1\sigma _y^2 \sigma _x + ~\mbox{}_1\sigma _y^3)$\par $~
\mbox{}_{\sigma }\int_0^{\infty } ~ \mbox{}_{\sigma }\int_0^{\infty
} \{ {\sf F}(\frac{(b,x+\eta +y
>}{2}) [{\sf F}(\frac{(b,2x+\zeta +\eta >}{2}) {\sf K}(x,x+\zeta
)]d\zeta \} d\eta $. \\ To simplify appearing formulas we use the
identity
\par $(28)$ $(\mbox{}^1L_2 + ~\mbox{}^2L_2) \{ {\sf F}(\frac{(b,x+\eta +y
>}{2}) [{\sf F}(\frac{(b,2x+\zeta +\eta >}{2}) g(x,\xi
)] \} =$\par $ \{ (\mbox{}^1\sigma _x + ~\mbox{}^2\sigma
_x)[\mbox{}^1\sigma _x ~\mbox{}^2\sigma _x + ~ \mbox{}^2\sigma _x
~\mbox{}^1\sigma _x] + 3\mbox{}^1_1\sigma _y^2~\mbox{}^2\sigma _x +
3\mbox{}^1_1\sigma _y (\mbox{}^1\sigma _x ~\mbox{}^2\sigma _x +
~\mbox{}^2\sigma _x ~\mbox{}^1\sigma _x + ~\mbox{}^2\sigma _x^2)
\}$\par $ \{ {\sf F}(\frac{(b,x+\eta +y
>}{2}) [{\sf F}(\frac{(b,2x+\zeta +\eta >}{2}) g(x,\xi
)] \} $\par $= 3\sigma _{\eta } (~\mbox{}^2\sigma _z
~\mbox{}^1\sigma _{\zeta } +  ~\mbox{}^1\sigma _{\zeta
}~\mbox{}^2\sigma _z) \{ {\sf F}(\zeta )|_{\zeta =\frac{(b,x+\eta +y
>}{2}} [{\sf F}(z)|_{z=\frac{(b,2x+\zeta +\eta >}{2}} g(x,\xi
)] \} $ \\  for every ${\cal A}_r$ valued function $g(x,\xi ),$
since the function ${\sf F}$ satisfies Conditions $(25)$ and \par
$(28.1)$ $\sigma _x ~\mbox{}_1\sigma _y^2 [f(x,y)g(x,y)] =
(\mbox{}^1\sigma _x + ~\mbox{}^2\sigma _x) [(\mbox{}^1_1\sigma _y
~\mbox{}^2_1\sigma _y + ~\mbox{}^2_1\sigma _y ~\mbox{}^1_1\sigma _y)
+ (\mbox{}^1_1\sigma _y^2 + ~\mbox{}^2_1\sigma _y^2)]
[f(x,y)g(x,y)]$ and \par $(28.2)$ $\sigma _x^3 [f(x,y)g(x,y)] =
(\mbox{}^1\sigma _x + ~\mbox{}^2\sigma _x) [(\mbox{}^1\sigma _x
~\mbox{}^2\sigma _x + ~\mbox{}^2\sigma _x ~\mbox{}^1\sigma _x) +
(\mbox{}^1\sigma _x^2 + ~\mbox{}^2\sigma _x^2)] [f(x,y)g(x,y)]$
\\ for each ${\cal A}_r$ valued twice differentiable
functions $f$ and $g$.
\par Then we have also from $(3,5)$:
\par $(29)$ $(\sigma _x-\sigma _y){\sf K}_2(x,y) = (~\mbox{}^1\sigma _x- ~\mbox{}^1\sigma _y
+ ~\mbox{}^2\sigma _x) ~\mbox{}_{\sigma }\int_0^{\infty } [{\sf
F}(\frac{(b,x+\zeta +y>}{2}) {\sf K}(x,x+\zeta) ] d\zeta $ \par $=
~\mbox{}^2\sigma _x ~\mbox{}_{\sigma }\int_0^{\infty } [{\sf
F}(\frac{(b,x+\zeta +y>}{2}) {\sf K}(x,x+\zeta) ] d\zeta $, \\ since
$\partial _{x_j}{\sf F}(x,y)=\partial _{y_j}{\sf F}(x,y)$ for each
$x, y$ and $j$.

\par Hence, the identity is satisfied:
\par $(30)$ $\sigma _x ~\mbox{}^2\sigma _x [~\mbox{}_{\sigma } \int_0^{\infty } {\sf
F}(\frac{(b, 2x+\zeta >}{2}) {\sf K}(x,x+\zeta ) d\zeta ] = - 2
~\mbox{}^2 \sigma _x [{\sf K}(x,x) {\sf K}(x,x)] - 2
(~\mbox{}^1\sigma _x- ~\mbox{}^1\sigma _y) [{\sf K}(x,y) {\sf
K}(x,x)]|_{y=x} - 2 [\sigma _x,\sigma _y] {\sf K}_2(x,y)|_{y=x}$, \\
since $[(\sigma _x+\sigma _y),(\sigma _x-\sigma _y)]= - 2 [\sigma
_x, \sigma _y]$ and
\par $\sigma _x ~\mbox{}^2\sigma _x [~\mbox{}_{\sigma } \int_0^{\infty } {\sf
F}(\frac{(b, 2x+\zeta >}{2}) {\sf K}(x,x+\zeta ) d\zeta ] = \sigma
_x\{ [(\sigma _x-\sigma _y){\sf K}_2(x,y)]|_{y=x} \} $ \par $ =
(\sigma _x+\sigma _y) (\sigma _x-\sigma _y) {\sf K}_2(x,y)|_{y=x}
$\par $=(\sigma _x-\sigma _y) (\sigma _x+\sigma _y) {\sf
K}_2(x,y)|_{y=x} - 2[\sigma _x,\sigma _y] {\sf K}_2(x,y)|_{y=x}$
\par $=(\sigma _x-\sigma _y)[-2{\sf K}(x,y){\sf K}(x,x)]|_{y=x} -
2[\sigma _x,\sigma _y] {\sf K}_2(x,y)|_{y=x}$ \\ due to Identity
$(16)$. On the other hand,
\par $[\sigma _x,\sigma _y] {\sf K}_2(x,y) =
(\mbox{}^1\sigma _x ~\mbox{}^1\sigma _y - ~\mbox{}^1\sigma _y
~\mbox{}^1\sigma _x + ~\mbox{}^2\sigma _x ~\mbox{}^1\sigma _y -
~\mbox{}^1\sigma _y ~\mbox{}^2\sigma _x) ~\mbox{}_{\sigma
}\int_0^{\infty } [{\sf F}(\frac{(b,x+\zeta +y>}{2}) {\sf
K}(x,x+\zeta) ] d\zeta $
\par $=(~\mbox{}^2\sigma _x ~\mbox{}^1\sigma _y -
~\mbox{}^1\sigma _y ~\mbox{}^2\sigma _x) ~\mbox{}_{\sigma
}\int_0^{\infty } [{\sf F}(\frac{(b,x+\zeta +y>}{2}) {\sf
K}(x,x+\zeta) ] d\zeta $, \\ since $(\mbox{}^1\sigma _x
~\mbox{}^1\sigma _y - ~\mbox{}^1\sigma _y ~\mbox{}^1\sigma _x) [{\sf
F}(\frac{(b,x+\zeta +y>}{2}) {\sf K}(x,x+\zeta) ]=0$, consequently,
\par $[\sigma _x,\sigma _y] {\sf K}_2(x,y)|_{y=x} =
\frac{1}{2} (~\mbox{}^2\sigma _x ~\mbox{}^1\sigma _x -
~\mbox{}^1\sigma _x ~\mbox{}^2\sigma _x) ~\mbox{}_{\sigma
}\int_0^{\infty } [{\sf F}(\frac{(b,2x+\zeta >}{2}) {\sf
K}(x,x+\zeta) ] d\zeta $ and \par $(31)$ $(\sigma _x
~\mbox{}^2\sigma _x + ~\mbox{}^2\sigma _x \sigma _x)
[~\mbox{}_{\sigma } \int_0^{\infty } {\sf F}(\frac{(b, 2x+\zeta
>}{2}) {\sf K}(x,x+\zeta ) d\zeta ] = - 4 ~\mbox{}^2 \sigma _x [{\sf
K}(x,x) {\sf K}(x,x)] - 4 (~\mbox{}^1\sigma _x- ~\mbox{}^1\sigma _y)
[{\sf K}(x,y) {\sf K}(x,x)]|_{y=x}$.
\par  Using Equations $(28-28.2)$ we rewrite Equation
$(27)$ in accordance with 4$(9-11)$ in the form:
\par $(32)$ $(I-A_x) [L_2{\sf K}(x,y)] =$\par $ \frac{\sf p}{4}
(L_2 - ~\mbox{}^3L_2) ~ \mbox{}_{\sigma }\int_0^{\infty } ~
\mbox{}_{\sigma }\int_0^{\infty } \{ {\sf F}(\frac{(b,x+y+\eta
>}{2}) [{\sf F}(\frac{(b,2x+\zeta +\eta >}{2}) {\sf K}(x,x+\zeta )]d\eta \}
d\zeta $ \par $+\frac{3\sf p}{4} [(\mbox{}^1\sigma _x ~
\mbox{}^3\sigma_x + ~ \mbox{}^3\sigma _x ~ \mbox{}^1\sigma _x+
~\mbox{}^2\sigma _x ~ \mbox{}^3\sigma_x + ~ \mbox{}^3\sigma _x ~
\mbox{}^2\sigma _x ) \mbox{}_1^1\sigma _y + ~ \mbox{}^1_1\sigma_y^2
~ \mbox{}^3\sigma _x + (\mbox{}^1\sigma _x ~ \mbox{}^2\sigma_x + ~
\mbox{}^2\sigma _x ~ \mbox{}^1\sigma _x+ ~\mbox{}^1\sigma _x^2 + ~
\mbox{}^2\sigma_x^2) ~ \mbox{}^3\sigma _x]$\par $ ( ~
\mbox{}_{\sigma }\int_0^{\infty } ~ \mbox{}_{\sigma }\int_0^{\infty
} \{ {\sf F}(\frac{(b,x+y+\eta >}{2}) [{\sf F}(\frac{(b,2x+\zeta
+\eta >}{2}) {\sf K}(x,x+\zeta )]d\eta \} d\zeta$
\par $+\frac{3\sf p}{4} [(\mbox{}^1\sigma _x + ~
\mbox{}^2\sigma_x + ~ \mbox{}^1\sigma _y) ~ \mbox{}^3\sigma _x^2 ( ~
\mbox{}_{\sigma }\int_0^{\infty } ~ \mbox{}_{\sigma }\int_0^{\infty
} \{ {\sf F}(\frac{(b,x+y+\eta>}{2}) [{\sf F}(\frac{(b,2x+\zeta
+\eta >}{2})$\par $ {\sf K}(x,x+\zeta )]d\eta \} d\zeta $
\par $= - \frac{3\sf p}{4} [\mbox{}^1\sigma _z  ~
\mbox{}^2\sigma_{\xi } + ~ \mbox{}^2\sigma _{\xi } ~\mbox{}^1\sigma
_z] ~ \mbox{}_{\sigma } \int_0^{\infty }  {\sf
F}(z)|_{z=\frac{(b,x+y
>}{2}} [{\sf F}(\xi )|_{\xi = \frac{(b,2x+\zeta>}{2}} {\sf K}(x,x+\zeta )]
d\zeta $
\par $ - \frac{3\sf p}{2} ~~\mbox{}^3\sigma _x^2 ~ \mbox{}_{\sigma } \int_0^{\infty }  {\sf
F}(\frac{(b, x+y>}{2}) [{\sf F}(\frac{(b,2x+\zeta>}{2}) {\sf
K}(x,x+\zeta )] d\zeta $ \par $+ 3{\sf p} ~~\mbox{}^3 \sigma _x
\sigma _{\eta }^2 ~ \mbox{}_{\sigma } \int_0^{\infty } ~
\mbox{}_{\sigma } \int_0^{\infty }  \{ {\sf F}(\frac{(b,x+\eta + y
>}{2}) [{\sf F}(\frac{(b,2x+\zeta +\eta>}{2}) {\sf K}(x,x+\zeta )]
d\eta \} d\zeta $
\\ since
\par $(\sigma _x^2+ ~\mbox{}_1\sigma _y^2) [{\sf F}
(\frac{(b, x+\eta +y>}{2}) {\sf F}(\frac{(b,2x+\zeta +\eta>}{2})] =
4\sigma _{\eta }^2 [{\sf F} (\frac{(b, x+\eta +y>}{2}) {\sf
F}(\frac{(b,2x+\zeta +\eta>}{2})]$.
\par If ${\sf F}\in Mat_s({\bf R})$ and ${\sf K}\in
Mat_s ({\cal A}_r)$, $2\le r\le 3$, $s\in {\bf N}$ for $r=2$, $s=1$
for $r=3$, then using Equation $(11)$ we deduce that \par $(33)$
$(I-A_x)[L_2{\sf K}(x,y)]= - \frac{3\sf p}{4} ( \mbox{}^1\sigma
_{\eta } ~\mbox{}^2\sigma _x + ~ \mbox{}^2\sigma _x ~\mbox{}^1\sigma
_{\eta })[{\sf F}(\eta )|_{\eta = \frac{x+y}{2}} {\sf K}_2(x,x)] -
\frac{3\sf p}{2} ~ \mbox{}^3\sigma _x (\mbox{}^2\sigma _x + ~
\mbox{}^3\sigma _x) \{ {\sf F}(\frac{(b,x+y>}{2}) ~ \mbox{}_{\sigma
} \int_0^{\infty }[{\sf F}(\frac{(b,2x+\zeta >}{2}) {\sf
K}(x,x+\zeta )]d\zeta \} $.
\par From Equations $(16,31,33)$   the identity
\par $(34)\quad (I-A_x)[L_2{\sf K}(x,y)]= \{ \sigma
_{\xi }  [ - \frac{3\sf p}{2} (I-A_x) (\sigma _x+~\mbox{}_1\sigma
_y) {\sf K}(x,y)$\par $ - \frac{3{\sf p}^2}{4} [(I-A_x) {\sf
K}(x,y)]{\sf K}_2(x,x) ] {\sf K}_2(\xi ,\xi ) \} |_{\xi =x}$
\par $- \frac{3\sf p}{2} \{ {\sf
F}(\frac{(b,x+y>}{2}) (-2~\mbox{}^2\sigma _x) ({\sf K}(x,x) {\sf
K}(x,x) ) - 2 [(\mbox{}^2\sigma _x - ~\mbox{}^2_1\sigma _y ) {\sf
F}(\frac{(b,x+y>}{2}) ({\sf K}(x,y) {\sf K}(x,x) )]|_{y=x} \} $ \\
follows. Expressing ${\sf F}$ from Equality $(1)$: \par ${\sf
F}(\frac{(b,x+y>}{2}) ={\sf F}(x,y) = (I-A_x){\sf K}(x,y)$ \\
and using Formulas $(21)$ and 4$(7,8)$ and the invertibility of the
operator $(I-A_x)$ and 7$(13)$ (see Remark 3.2 also) we get
\par $(35)$ $L_2{\sf K}(x,y) =
\{ [ - \frac{3\sf p}{2} (\sigma _x+ ~\mbox{}_1\sigma _y) {\sf
K}(x,y) - \frac{3{\sf p}^2}{4} {\sf K}(x,y){\sf K}_2(x,x) ] (- {\sf
K}^2(x,x) \} $
\par $- \frac{3\sf p}{2} \{  {\sf
K}(x,y)[ (-2~ \mbox{}^2\sigma _x) ({\sf K}(x,x) {\sf K}(x,x) )] -
2{\sf K}(x,y)[(\mbox{}^1\sigma _x - ~\mbox{}^1_1\sigma _y )  ({\sf
K}(x,y) {\sf K}(x,x) )]|_{y=x} \} $
\par $=3{\sf p} (\mbox{}^1\sigma _x+~\mbox{}_1\sigma _y) {\sf K}(x,y){\sf
K}^2(x,x) +\frac{3{\sf p}^2}{2} [{\sf K}(x,y){\sf K}_2(x,x)]{\sf
K}^2(x,x)$\par $+ 3{\sf p} {\sf K}(x,y) [~\mbox{}^2\sigma _x ({\sf
K}(x,x){\sf K}(x,x))] +3{\sf p} {\sf K}(x,y) [ (\mbox{}^1\sigma _x -
~ \mbox{}^1_1\sigma _y) ({\sf K}(x,y) {\sf K}(x,x))]|_{y=x}$.

\par Using Formula $(22)$ permits to simplify Formula $(35)$ in the
considered case to: \par $(36)$ $L_2{\sf K}(x,y)=3{\sf p} [(\sigma
_x+ \sigma _y) {\sf K}(x,y)]{\sf K}^2(x,x) +3{\sf p} {\sf
K}(x,y)[\mbox{}^2\sigma _x ({\sf K}(x,x) {\sf K}(x,x))]$.
\par
Putting $g(x,t) ={\sf K}(x,x,t)$ one gets particularly for $\sigma
=\mbox{}_1\sigma $ from Equation $(36)$:
\par $(37)$ $(\mbox{}_2\sigma _t+\sigma _x^3)g(x,t) = 3{\sf p} [\sigma
_xg(x,t)]g^2(x,t) + 3{\sf p} g(x,t)[\mbox{}^2\sigma
_x(g(x,t)g(x,t))]$.
\par The latter two partial differential equations are non-linear
even over the real field. Particularly, if $\mbox{}_2\sigma
_t=\partial /\partial t_0$ and $\sigma _x=\partial /\partial x_0$
Equation $(37)$ gives the modified Korteweg-de-Vries equation:
\par $\partial g(x,t)/\partial t +\partial ^3g(x,t)/\partial x^3- 6pg^2(x,t)
\partial g(x,t)/\partial x=0$ with $t=t_0$ and $x=x_0$.
\par Find particular solutions of $(36)$, when \par $\mbox{}_2\sigma
_tf=(\partial f/\partial t_0)$ and $\sigma _x g(x)=
\sum_{n_0+1}^{n+n_0}i_j^*[\partial g(x)/\partial x_j]$, using
\par ${\sf F}(\frac{(b,x+y>}{2}),t) := (nt)^{-1/3} f(\xi _{n_0+1}+\eta
_{n_0+1},...,\xi _{n+n_0}+\eta _{n+n_0})$ with $\xi = x(nt)^{-1/3}$
and $\eta =y(nt)^{-1/3}$, $~ -1\le n_0$, $~2\le n$, $~n_0+n\le 3$,
$t=t_0$. While with the help of the operators
\par $L_{2,j}g(x,y,t)=\frac{1}{n}i_j^* \partial _tg(x,y,t)+(\partial
_{x_j}\sigma _x^2 + \partial _{y_j}\sigma _y^2 + 3\partial
_{x_j}\sigma _y^2+3\partial _{y_j}\sigma _x^2)g(x,y,t)$ and the
conditions $L_{2,j}F=0$ for each $j=n_0+1,...,n_0+n$ one obtains the
partial differential equations
\par $f +\xi _j f +\partial _{\xi _j}\sum_{n_0+1}^{n+n_0} \partial ^2_{\xi
_k}f=0$ \\ for each $j=n_0+1,....,n+n_0$, where $f=f(\xi
_{n_0+1},...,\xi _{n_0+n})$, $\partial _{x_j}=\partial /\partial
_{x_j}$. Solving these partial differential equations we get the
solution ${\sf K}(x,y) = (nt)^{-1/3} \kappa (\xi _{n_0+1},...,\xi
_{n_0+n},\eta _0,...,\eta _{n_0+n})$ from the integral equation
\par ${\hat K}(\xi ,\eta ) = {\hat F}(\xi +\eta ) + \frac{p}{4} ~
(~ \mbox{}_{\sigma }\int_{\xi }^{\infty }(~ \mbox{}_{\sigma
}\int_{\xi }^{\infty } {\hat F}(u+\eta )[{\hat F}(z+u){\hat K}(\xi ,z)]dz)du$, \\
where ${\hat K}(\xi ,\eta )=\kappa (\xi _{n_0+1},...,\xi
_{n_0+n},\eta _0,...,\eta _{n_0+n})$, $\xi =\xi _0i_0+...+\xi
_3i_3$.
\par Each solution of the modified Korteweg-de-Vries equation
generates a solution of Kadomtzev-Petviashvili partial differential
equation due to Miura's transform: $g=-v^2-\sigma _xv$. In the
considered non-commutative case it takes the form:
\par $\mbox{}_2\sigma _t(v^2) + ~ \mbox{}_2\sigma _t\sigma _x v +
(\sigma _x^3(v^2)+\sigma _x^4v) - 3{\sf p} (\sigma _x(v^2) +\sigma
_x^2v) (v^4+v^2\sigma _xv + (\sigma _xv) v^2+(\sigma _xv)^2) -3{\sf
p} (v^2+ \sigma _xv)[\mbox{}^2\sigma _x ((v^2+\sigma _xv)(v^2+\sigma
_xv)) ]$.
\par Thus we have demonstrated the following statement.
\par {\bf 9.1 Theorem.} {\it Partial differential Equation $(36)$
with $\psi _0=~\mbox{}_1\psi _0=0$ over the Cayley-Dickson algebra
${\cal A}_r$, $2\le r\le 3$, has a solution given by Formulas
$(1-3,25,26)$, when the appearing integrals uniformly converge by
parameters as in Proposition 4 and the operator $(I-A_x)$ is
invertible and ${\sf F}\in Mat_s({\bf R})$ and ${\sf K}\in
Mat_s({\cal A}_r)$ with $s\in {\bf N}$ for $r=2$ and $s=1$ for
$r=3$.}

\par {\bf 10. Example.} Now we consider the pair of partial
differential operators \par $(1)$ $L_1=\sigma _x+ \sigma _y$ and
\par $(2)$ $L_2= \mbox{}_1\sigma _t +\sigma _x^2 + u \sigma _y\sigma _x
+\sigma _y^2$, \\ where $u\in {\bf R}$ is a real constant, and the
integral equation:
\par $(3)$ ${\sf K}(x,y) = {\sf F}(x,y) + {\sf p} ~\mbox{}_{\sigma
}\int_x^{\infty } {\sf F}(z,y) {\sf K}(x,z)dz$. \par Acting on the
both sides of the latter equation with the operator $L_2$ and using
the condition \par $(5)$ $L_{2,j}{\sf F}(x,y)=0$ for each $j$, where
\par $(6)$ $L_{2,j}{\sf F}(x,y)= [\mbox{}_1\psi _j\partial _{t_j} + $\par $u
\sum_{s,p: ~i_si_p=i_j} \psi _s\psi _p(\partial _{y_s} \partial
_{x_p} + (i_pi_s)(i_si_p) \partial _{y_p}\partial _{x_s})+ \pi _j
\circ (\sigma _x^2 +\sigma _y^2)] {\sf F}(x,y)$, \\ where $\pi _j:
X\to X_j$ is the ${\bf R}$ linear projection operator (see \S 1.2).
We infer that
\par $(7)$ $L_2{\sf K}(x,y) ={\sf p} (\mbox{}_1^1\sigma _t +\sigma _x^2 + u \sigma _y\sigma _x
+\sigma _y^2 + \mbox{}_1^2\sigma _t)~\mbox{}_{\sigma }\int_x^{\infty
} {\sf F}(z,y) {\sf K}(x,z)dz$
\par $= {\sf p} (\mbox{}_1^2\sigma _t +\sigma _x^2 + u ~\mbox{}^1\sigma _y\sigma _x
- ~\mbox{}^1\sigma _z^2 - u ~\mbox{}^1\sigma _y~\mbox{}^1\sigma _z)
~\mbox{}_{\sigma }\int_x^{\infty } {\sf F}(z,y) {\sf K}(x,z)dz$, \\
since $(\mbox{}_1\sigma _t +\sigma _y^2) {\sf F}(z,y) = - (\sigma
_z^2 + u \sigma _y\sigma _z){\sf F}(z,y)$. In view of Proposition 4
and Corollary 5 we get the formula:
\par \par $(8)$ $L_2{\sf K}(x,y) = I_1+I_2$, \\ where
\par $(9)$ $I_1= {\sf p} (~\mbox{}_1^2\sigma _t+ \sigma _x^2 + u~\mbox{}^1\sigma
_y~\sigma _x)~\mbox{}_{\sigma }\int_x^{\infty } {\sf F}(z,y) {\sf
K}(x,z)dz$\par $ = {\sf p} ( \mbox{}^2\sigma _x^2 +
u~\mbox{}^1\sigma _y~\mbox{}^2\sigma _x) ~\mbox{}_{\sigma
}\int_x^{\infty } {\sf F}(z,y) {\sf K}(x,z)dz $\par $ - {\sf p}
\sigma _x[{\sf F}(x,y){\sf K}(x,x)] - {\sf p} ~\mbox{}^2 \sigma _x
[{\sf F}(x,y){\sf K}(x,z)]|_{z=x} - u {\sf p} \sigma _y [{\sf
F}(x,y){\sf K}(x,x)]$ and
\par $(10)$ $ ~ -I_2 = {\sf p} (\mbox{}^1\sigma _z^2 +
u~\mbox{}^1\sigma _y~\mbox{}^1\sigma _z) ~\mbox{}_{\sigma
}\int_x^{\infty } {\sf F}(z,y) {\sf K}(x,z)dz$
\par $ = {\sf p} (\mbox{}^2\sigma _z^2 - u ~\mbox{}^1\sigma _y~\mbox{}^2\sigma _z)
~\mbox{}_{\sigma }\int_x^{\infty } {\sf F}(z,y) {\sf K}(x,z)dz $\par
$ - {\sf p} ~\mbox{}^1\sigma _x[{\sf F}(x,y){\sf K}(x,x)] +{\sf p}
~\mbox{}^2 \sigma _z [{\sf F}(x,y){\sf K}(x,z)]|_{z=x} - u{\sf p}
\sigma _y [{\sf F}(x,y){\sf K}(x,x)]$. \par Then using the second
condition
\par $(11)$ $L_1{\sf F}(x,y)=0$, i.e. ${\sf F}(x,y) = {\sf
F}(\frac{(\psi , x-y>}{2})$, \\ it is possible to simplify certain
terms appeared above:
\par $(12)$ $\sigma _y ~\mbox{}^2\sigma _x ~\mbox{}_{\sigma
}\int_x^{\infty } {\sf F}(z,y) {\sf K}(x,z)dz $
\par $=  ~\mbox{}^2\sigma _z ~\mbox{}^2\sigma _x ~\mbox{}_{\sigma
}\int_x^{\infty } {\sf F}(z,y) {\sf K}(x,z)dz + ~\mbox{}^2\sigma _x
[{\sf F}(x,y) {\sf K}(x,z)]|_{z=x}$ and
\par $(13)$ $\mbox{}^1\sigma _y ~\mbox{}^2\sigma _z ~\mbox{}_{\sigma
}\int_x^{\infty } {\sf F}(z,y) {\sf K}(x,z)dz = - ~\mbox{}^1\sigma
_z ~\mbox{}^2\sigma _z ~\mbox{}_{\sigma }\int_x^{\infty } {\sf
F}(z,y) {\sf K}(x,z)dz $
\par $=   ~\mbox{}^2\sigma _z^2 ~\mbox{}_{\sigma
}\int_x^{\infty } {\sf F}(z,y) {\sf K}(x,z)dz
 + ~\mbox{}^2\sigma _z [{\sf F}(x,y) {\sf K}(x,z)]|_{z=x}$.
\par Therefore, Equations $(7-13)$ lead to the identity:
\par $(14)$ $L_2{\sf K}(x,y) = {\sf p} (\mbox{}^2_1\sigma _t+ ~\mbox{}^2\sigma _x^2
+ u ~~\mbox{}^2\sigma _z ~\mbox{}^2\sigma _x  + (u-1)
~~\mbox{}^2\sigma _z^2 ) ~\mbox{}_{\sigma }\int_x^{\infty } {\sf
F}(z,y) {\sf K}(x,z)dz $
\par $- (u+2){\sf p} ~~\mbox{}^2\sigma _x [{\sf F}(x,y) {\sf K}(x,x)] +u{\sf p}
~~\mbox{}^2\sigma _x[{\sf F}(x,y) {\sf K}(x,z)]|_{z=x} +u{\sf p}~
~\mbox{}^2\sigma _z[{\sf F}(x,y) {\sf K}(x,z)]|_{z=x}$.
\par Take $u=2$. If ${\sf F}(x,y)\in
Mat_s({\bf R})$ and ${\sf K}(x,y)\in Mat_s({\cal A}_r)$ for each
$x,y \in U$, $2\le r \le 3$, with $s\in {\bf N}$ for $r=2$ and $s=1$
for $r=3$, then Formulas $(3,14)$ and 2$(7)$ and 7$(2,13)$ and
5$(7)$ for the invertible operator $(I-A_x)$ imply, that
\par $(15)$ $(\mbox{}_1\sigma _t + \sigma _x^2
+ 2 ~\mbox{}^2\sigma _y ~\mbox{}^2\sigma _x  + \sigma _y^2 ) {\sf
K}(x,y) = - 2 {\sf p}{\sf K}(x,y) [\sigma _x {\sf K}(x,x)] $.
\par Putting $g(x,t)={\sf K}(x,x)$ we get particularly on the
diagonal $x=y$:
\par $(16)$ $(\mbox{}_1\sigma _t + \sigma _x^2) g(x,t) =
- 2 {\sf p}  g(x,t)[\sigma _x g(x,t)]$.
\par If
$u\in {\cal I}_r: = \{ z\in {\cal A}_r: ~ Re(z)=0 \} $ and $\psi
_0=0$, then $\mbox{}^2\sigma _x (u(x,t) u(x,t)) + u(x,t) (\sigma _x
u(x,t))= \sum_{j>0} u_j(x,t)(\partial u(x,t)/\partial x_j) \psi _j$.
\par This approach has the following application.
\par A non-isothermal flow of a non-compressible
Newtonian liquid with a dissipative heating is described by the
system of partial differential equations:
\par $(17)$ $\rho \sum_{j=1}^3 u_j \partial u/\partial x_j =
- \frac{1}{2} \mu \sigma _x^2 u +\rho \frac{\partial u}{\partial
t_0}$
\par $(18)$ $div (u) = \sum_{j=1}^3 \partial u_j/\partial x_j=0$
\par $(19)$ $\rho c_p\sum_{j=1}^3 u_j
\partial T/\partial x_j = - \lambda \sigma _x^2 T + 2\mu I_2$, \\
when a density $\rho $ and a dynamical viscosity $\mu $ are
independent of coordinates $x_1, x_2, x_3$, where
\par $I_2=\frac{1}{4} \sum_{j,k=1}^3 (\partial u_k/\partial x_j +
\partial u_j/\partial x_k)^2$ \\ denotes the tensor of deformation
velocities, $T$ is the temperature, also $\lambda $ and $c_p$ are
physical constants. Here we take $u=u_1i_1+u_2i_2+u_3i_3$ with
real-valued functions $u_1, u_2, u_3$ and $\sigma
_xf(x)=\sum_{j=1}^3i_j^* (\partial f(x)/\partial x_j)$,
$~x=x_0i_0+x_1i_1+x_2i_2+...+x_{2^r-1}i_{2^r-1}\in {\cal A}_r$,
$r=2$ or $r=3$, $~x_0$ and $x_j$ with $j>3$ can be taken constant,
for example, $x_0=0$, $x_j\in {\bf R}$ for each $j$,
$~\mbox{}_1\sigma _t =\mbox{}_1\psi _0 ~
\partial_{t_0}$ and choose suitable constants $\mbox{}_1\psi _0$ and ${\sf p}$.
\par For usual liquids like water the dependence of the dynamical
viscosity on the temperature is described as $\mu =\mu _0(T_0/T)^m$,
where $m\ge 0$. For very viscous Newtonian liquids like glycerin the
function $\mu =\mu _0 \exp [-\beta (T-T_0)]$ is usually taken with
empirical constants $\mu _0, T_0$ and $\beta $ (see also
\cite{polzayrb}). Using the method described above it is possible to
find solutions satisfying Equation $(20)$ and take into account
Condition $(21)$ with the help of ${\sf F}(x,x)$ and then ${\sf
K}(x,x)$ satisfying it. \par Put $g=v+u$, where $u\in {\cal I}_2 :=
\{ z \in {\cal A}_2: ~ Re (z)=0 \} $ and $v\in {\cal A}_r\ominus
{\cal I}_2$, where $2\le r \le 3$. If $g$ is a solution of partial
differential equation $(16)$ and \par $v(\sigma u)+u(\sigma v) =
~\mbox{}^2\sigma (uu + vv)$, \\ then $u$ is a solution of partial
differential equation $(17)$. Condition $(18)$ means that $Re
(\sigma u)=0$ or $\sigma u+ (\sigma u)^*=0$. Calculating $I_2$ and
solving Equation $(19)$ one also calculates the temperature.
\par The result of this section can be formulates as the following
theorem.
\par {\bf 10.1 Theorem.} {\it Partial differential Equation $(15)$
over the Cayley-Dickson algebra ${\cal A}_r$ with $2\le r\le 3$ has
a solution given by Formulas $(1,3,5,6,11)$, when the appearing
integrals uniformly converge by parameters on each compact
sub-domain (see Proposition 4) and the operator $(I-A_x)$ is
invertible and ${\sf F}\in Mat_s({\bf R})$ and ${\sf K}\in
Mat_s({\cal A}_r)$, $s\in {\bf N}$ for $r=2$, $s=1$ for $r=3$.}

\par {\bf 11. Remark.} The results presented above show that the method of
commutators of integral operators with partial differential
operators over complex numbers becomes more powerful with the use of
Cayley-Dickson algebras and the non-commutative line integral over
them and permits to solve more general non-linear  partial
differential equations. It is planned to develop further this method
for partial differential equations with variable coefficients and
also for generalized functions and equations with generalized
coefficients.

\par Department of Applied Mathematics,
\par Moscow State Technical University MIREA,
av. Vernadsky 78,
\par Moscow, Russia
\par e-mail: sludkowski@mail.ru
\end{document}